\documentclass[paper=letter, fontsize=12pt]{scrartcl}
\usepackage{fourier,amssymb,epsfig,psfrag,amsmath,art,wrapfig,multicol,color,epstopdf,float,cases}
\usepackage{subfigure}[1995/03/06]

\def\IR{\mathbb{R}}
\newcommand{\rom}{\mathrm}
\newcommand{\mT}{^\mathrm{T}}
\def\norm#1{\|#1\|}

\newtheorem{definition}{Definition}
\newtheorem{assumption}{Assumption}
\newtheorem{lemma}{Lemma}
\newtheorem{theorem}{Theorem}

\newtheorem{remark}{Remark}

\newcommand{\horrule}[1]{\rule{\linewidth}{#1}}

\DeclareMathAlphabet{\mathcal}{OMS}{cmsy}{m}{n}
\usepackage[english]{babel}
\begin{document}
\definecolor{tBlue}{RGB}{0,0,0}
\newcommand{\blue}[1]{\color{tBlue}{#1}}
\newcommand{\phid}{\phi_\mathrm{d}\left(\left\|e_i(t)\right\|_{P_i}\right)}
\newcommand{\phidx}{\phi_\mathrm{d}\left(\left\|x_i(t)\right\|_{P_i}\right)}
\newcommand{\phie}{\phi\left(\left\|e_i(t)\right\|_{P_i}\right)}
\newcommand{\pePB}{\phid e_i^{\mathrm{T}}(t) P_iB_i}


\title{
\vspace{-0.6in}
\usefont{OT1}{bch}{b}{n}
\normalfont \normalsize \textsc{} \\ [25pt]
\horrule{2pt} \\[0.3cm] 
\Large \textbf{A Distributed Control Approach for  Heterogeneous \\ Linear Multiagent Systems}${^\star}$
{\thanks{\hspace{-0.15cm}$^\star$\hspace{0.08cm} This research was supported by the Dynamics, Control, and Systems Diagnostics Program of the National Science Foundation under Grant CMMI-1657637.}
\thanks{\hspace{-0.15cm}$^\ddagger$\hspace{0.08cm} Corresponding Author: Engineering Building C 2209, 4202 East Fowler Avenue, Tampa, Florida 33620, United States of America (Address); +1 813 974 5656 (Phone); yucelen@usf.edu (Email).}
} 
\horrule{2pt}
}
\def\GruenwaldYGDM
{\begin{tabular}{c} 
\textbf{Selahattin Burak Sarsilmaz}\hspace{0.05cm} and\hspace{0.05cm} \textbf{Tansel Yucelen}$^{\dagger}$\\
\textit{Department of Mechanical Engineering}\\
\textit{University of South Florida}\\
 \end{tabular}
} 
\author{\GruenwaldYGDM} \date{} \maketitle \baselineskip 16pt


\vspace{9.2cm}

\textbf{\textit{Abstract}} --- This paper considers an internal model based distributed control approach to the cooperative output regulation problem of heterogeneous linear time-invariant multiagent systems over fixed directed communication graph topologies. First, a new definition of the linear cooperative output regulation problem is introduced in order to allow a broad class of functions to be tracked and rejected by a network of agents. Second, the solvability of this problem with three distributed control laws, namely dynamic state feedback, dynamic output feedback with local measurement, and dynamic output feedback, is investigated by first considering a global condition and then providing an agent-wise local sufficient condition under standard assumptions. Finally, two numerical examples are provided to illustrate the selected contributions of this paper.

\baselineskip 16pt


\vspace{0.0em}

\textbf{\textit{Keywords}} --- Heterogeneous multiagent systems;\hspace{0.05cm} distributed control;\hspace{0.05cm} cooperative output regulation;\hspace{0.05cm} internal model

\clearpage 
\baselineskip=21.25pt 
\setcounter{page}{1}


\section{Introduction}\label{Introduction}

Heterogeneous multiagent systems formed by networks of agents having different dynamics and dimensions present a significantly broader class of multiagent systems than their heterogeneous and homogeneous counterparts that consist of networks of agents having different dynamics with the same dimension and identical dynamics, respectively.
Therefore, analysis and synthesis of distributed control approaches for this class of multiagent systems that rely on local information exchange has been an attractive research topic in the systems and control field over the last decade. 

In particular, the cooperative output regulation problem of heterogeneous (in dynamics and dimension) linear time-invariant multiagent systems, where the output of all agents synchronize to the output of the leader, over general fixed directed communication graph topologies have been recently investigated in  \cite{Su_2012,su_huang_2012_b,Huang_C_tac,Li_2015,Farnaz_automatica,Cai_2017, Lu_Liu_TAC_2017}. This problem can be regarded as the generalization of the linear output regulation problem given in, for example, \cite{Huang_book} to multiagent systems. 
As a consequence, distributed control approaches to this regulation problem can be classified into two categories: 
\begin{itemize}
	\item The first category is predicated on feedforward design methodology, where the authors of \cite{Su_2012,su_huang_2012_b,Li_2015,Cai_2017, Lu_Liu_TAC_2017} present contributions. 
	In the presence of plant uncertainties, however, this methodology is known to be not robust
	since the feedforward gain of each agent relies on the solution of the regulator equations. 
	\item The second category is predicated on internal model principle, where the authors of \cite{Huang_C_tac,Farnaz_automatica} present contributions. 
	While this methodology is robust with respect to small variations of the plant parameters as compared to feedforward design methodology, it cannot be applied when the transmission zero condition does not hold.
\end{itemize}

The common denominator of
 these papers 
 is that an exosystem, which has an unforced linear time-invariant dynamics, generates both a reference trajectory and external disturbances to be tracked and rejected by networks of agents. 
Specifically, the system matrix of the exosystem is explicitly used by controllers of all agents in \cite{Su_2012,su_huang_2012_b,Li_2015, Lu_Liu_TAC_2017} and a proper subset of agents in \cite{Cai_2017}; or each agent incorporates a  $p$-copy internal model of this matrix in its controller \cite{Huang_C_tac,Farnaz_automatica}. 

\subsection{Contributions}\label{CONTR}

Considering applications of the distributed control approaches in \cite{Su_2012,su_huang_2012_b,Huang_C_tac,Li_2015,Farnaz_automatica,Cai_2017, Lu_Liu_TAC_2017}, it can be a challenge to precisely know the system matrix of the exosystem, even the dynamical structure of the exosystem; especially, when an external leader interacts with the network of agents or a control designer simply injects optimized trajectory commands to the network based on, for example, an online path planning algorithm. 
In order to guarantee ultimately bounded tracking error in such cases, a new, generalized definition for the cooperative output regulation problem is needed.

This paper focuses on heterogeneous (in dynamics and dimension) linear time-invariant multiagent systems over general fixed directed communication graph topologies. 
First, we present the generalized definition for the linear cooperative output regulation problem. 
Second, we investigate the solvability of this problem for internal model based distributed dynamic state feedback, output feedback with local measurement, and output feedback control laws. 
To this end, we not only consider global conditions but also provide agent-wise local sufficient conditions under standard assumptions. Considering large-scale applications of multiagent systems, the agent-wise local sufficient conditions are primarily important for
independent controller design of each agent (i.e., without depending on the dynamics of other agents ). 

The system-theoretical approach presented in this paper\footnote{\hspace{0.16cm}Although they are not completely related, \cite{sarsilmaz_yucelen_asme_dscc_2017,sarsilmaz_yucelen_acc_2018} may be regarded as preliminary works of this paper.} is relevant to the studies in \cite{Huang_C_tac,Farnaz_automatica}, where they also focus on the linear cooperative output regulation problem with an internal model based distributed dynamic state feedback control law. 
Specifically, \cite{Farnaz_automatica} extends the approach in \cite{Huang_C_tac} to an output feedback control under an output feedback stabilizability condition. 
In addition to the generalized definition of the linear cooperative output regulation problem, the contribution of this paper differs from the studies in \cite{Huang_C_tac,Farnaz_automatica} based on the following points: 
\begin{itemize}
	\item First, we note that the theoretical contribution of this paper covers not only the dynamic state feedback problem but also the dynamic output feedback problem with local measurement as well as the dynamic output feedback problem. 
	Unlike the results presented in \cite{Farnaz_automatica}, this paper does not assume the output feedback stabilizability for the dynamic output feedback problem with local measurement. 
	With regard to the dynamic output feedback problem, the results of this paper does not require agents to access their own states or outputs.  
	\item To prove the existence of a unique solution to the matrix equations that are crucial for the solvability of the problem, Section III in \cite{Huang_C_tac} (Theorem 4 in \cite{Farnaz_automatica}) decomposes these matrix equations, which consist of the overall dynamics of the multiagent system, into matrix equations, which deal with the dynamics of each agent separately.
	In contrast, we do not decompose these matrix equations; see the sixth paragraph of Appendix A for the advantage. 
	In particular, Lemma \ref{Lemma3} of this paper, which is also applicable to dynamic output feedback cases, guarantees that these matrix equations have a unique solution without requiring their decompositions. 
	\item A considerable number of gaps in the related results of  \cite{Huang_C_tac,Farnaz_automatica} is illustrated by counterexamples in Appendices and fixed in Appendices as well as in Section \ref{first_controller}.
\end{itemize}

\subsection{Organization}\label{Org}

The rest of the paper is organized as follows. Section \ref{P} presents the notation and the essential mathematical preliminaries. Section \ref{PS} formulates the linear cooperative output regulation problem considered in this paper. The solvability of this problem is investigated in Section \ref{DSF} and two illustrative numerical examples are presented in Section \ref{INE}. Finally, Section \ref{con} concludes the paper. 

\section{Mathematical Preliminaries}\label{P}

A standard notation is used in this paper. Specifically, $\IR$, $\IR^n$, and $\IR^{n \times m}$ respectively denote the sets of all real numbers, $n \times 1$ real column vectors, and $n \times m$ real matrices\footnote{\hspace{0.16cm}In this paper, all real matrices are defined over the field of complex numbers.}; $\textbf{1}_n$ and $I_n$ respectively denote the $n \times 1$ vector of all ones and the  $n \times n$ identity matrix; and ``$\triangleq$'' denotes equality by definition. 
We write $(\cdot)\mT$ for the transpose and $\norm{\cdot}_2$ for the induced two norm of a matrix; $\sigma(\cdot)$ for the spectrum\footnote{\hspace{0.16cm}We follow Definition 4.4.4 in \cite{Bernstein_book}.} and $\rho (\cdot)$ for the spectral radius of a square matrix; $(\cdot)^{-1}$ for the inverse of a nonsingular matrix; and $\otimes$ for the Kronecker product.
 We also write $A \leq B$ for $A \in \IR^{n \times m}$, $B \in \IR^{n \times m}$ if entries $a_{ij} 
	\leq b_{ij}$ for all ordered pairs $(i,j)$. 
 Finally,  $\text{diag}(A_1, \ldots, A_n)$ is a block-diagonal matrix with matrix entries $A_1, \dots, A_n$ on its diagonal.

We now concisely state the graph theoretical notation used in this paper, which is based on \cite{Lewis_book}. 
In particular, consider a fixed (i.e., time-invariant) directed graph ${\mathcal{G}} = (\mathcal{V} ,\mathcal{E})$, where $\mathcal{V} = \big\{v_1, \ldots, v_N \big\}$ is a nonempty finite set of $N$ nodes and $\mathcal{E} \subset \mathcal{V} \times \mathcal{V} $ is a set of edges. 
Each node in $\mathcal{V}$ corresponds to a follower agent.
There is an edge rooted at node $v_j$ and ended at $v_i$ (i.e., $(v_j, v_i) \in \mathcal{E} $)  if and only if $v_i$ receives information from $v_j$. $\mathcal{A} = [a_{ij}] \in \IR^{N \times N} $ denotes the adjacency matrix, which describes the graph structure; that is, $a_{ij} > 0 \Leftrightarrow (v_j, v_i) \in \mathcal{E}$ and $a_{ij} = 0$ otherwise. 
Repeated edges and self loops are not allowed; that is, $a_{ii} = 0, \ \forall i \in {\mathcal{N}}$ with  ${\mathcal{N}} = \big\{1,\ldots, N\big\}$. 
The set of neighbors of node $v_i$ is denoted as $N_i = \big\{j  \ | \ (v_j, v_i) \in \mathcal{E} \big\}$. 
In-degree matrix is defined by $\mathcal{D} = \rom{diag}(d_1, \ldots, d_N)$ with $d_i = \sum_{j\in N_i} a_{ij} $.
A directed path from node $v_i$ to node $v_j$ is a sequence of successive edges in the form $\big\{(v_i,v_p),(v_p,v_q), \ldots, (v_r,v_j)\big\}$. If $v_i = v_j$, then the directed path is called a loop.
A directed graph is said to have a spanning tree if there is a root node such that it has directed paths to all other nodes in the graph.
A fixed augmented directed graph is defined as  ${ \mathcal{\bar G}} = (\mathcal{\bar V} ,\mathcal{\bar E})$, where $ \mathcal{\bar V} = \big\{v_0,v_1, \ldots, v_N \big\}$ is the set of $N+1$ nodes, including leader node $v_0$ and all nodes in $\mathcal{V}$, and $\mathcal{\bar E} = \mathcal{E} \cup \mathcal{E'} $ is the set of edges with $\mathcal{E'}$ consisting of some edges in the form of $(v_0,v_i)$, $i \in \mathcal{N} $.

The concept of internal model introduced next slightly modifies Definition 1.22 and Remark 1.24 in \cite{Huang_book}.

\definition\label{Def1} Given any square matrix $A_0$, a triple of matrices $(M_1,M_2,M_3)$ is said to incorporate a $p$-copy internal model of the matrix $A_0$ if 
\begin{eqnarray}
M_1 = T  \begin{bmatrix}
S_1  & S_2 \\
0  & G_1  
\end{bmatrix}  T^{-1}, \ 
M_2 = T  \begin{bmatrix}
S_3 \\
G_2  
\end{bmatrix}, \  
M_3 =  T \begin{bmatrix}
S_4 \\
0  
\end{bmatrix},  \label{M1_M2_M3_Definition1_1} 
\end{eqnarray}
or 
\begin{eqnarray}
M_1 = G_1, \ M_2 = G_2, \ M_3 = 0, \label{M1_M2_M3_Definition1_2} 
\end{eqnarray}
where $S_l, \ \hspace{-0.06cm} l = 1, 2, 3, 4$, is any matrix with an appropriate dimension, $T$ is any nonsingular matrix with an appropriate dimension, the zero matrix in $M_3$ has as many rows as those of $G_1$, and 
\begin{eqnarray}
G_1 = \rom{diag}(\beta_1, \dots, \beta_p), \quad G_2 = \rom{diag}(\sigma_1, \dots, \sigma_p), \label{G1_G2_def} \nonumber
\end{eqnarray}	
 where for $l=1, \dots, p$, $\beta_l \in \IR^{s_l \times s_l}$ and $\sigma_l \in R^{s_l} $ satisfy the following conditions:

$a)$ The pair $(\beta_l,\sigma_l)$ is controllable. \\
\indent $b)$\hspace{0.0054cm} The minimal polynomial of $A_0$ is equal to the characteristic polynomial of $\beta_l$.

\section{Problem Formulation}\label{PS}

Consider a system of $N$ (follower) agents with heterogeneous linear time-invariant dynamics subject to external disturbances over a fixed directed communication graph topology $\mathcal{G}$. The dynamics of agent $i \in \mathcal{N}$ is given by
\begin{eqnarray}
\dot x_i(t) &=& A_ix_i(t)+B_iu_i(t) + \delta_i(t),    \quad x_i(0)=x_{i0},   \quad t  \geq 0, \label{local_state_1} \nonumber  \\ 
y_i(t)&=& C_ix_i(t) + D_iu_i(t), \label{local_output_1} \nonumber 
\end{eqnarray} 
with state $x_i(t)\in\IR^{n_i}$, input  $u_i(t)\in\IR^{m_i}$, output $y_i(t)\in\IR^{p}$, and external disturbance $\delta_i(t) = E_{\delta_i} \delta(t)\in\IR^{n_i}$, where $\delta(t) \in \IR^{q_{\delta}} $ is a solution to the unknown disturbance dynamics with an initial condition. In addition, the reference trajectory to be tracked is denoted by $y_0(t) = R_{\rom{r}} r_0(t) \in\IR^{p} $,  where $r_0(t) \in\IR^{q_\rom{r}} $ is a solution to the unknown leader dynamics with an initial condition.

Let $\omega(t) \triangleq [r_0^\rom{T}(t), \delta^\rom{T}(t)]\mT$ $\in\IR^{q}$ be the solution  of the unknown exosystem, where $q=q_\rom{r}+q_{\delta}$. Instead of assuming that the exosystem has an unforced linear time-invariant dynamics with a known system matrix (e.g., see \cite{Su_2012,Huang_C_tac,Farnaz_automatica}), we consider that the exosystem has an unknown dynamics.
From this perspective, the 
exosystem can represent any (e.g., linear or nonlinear) dynamics provided that its solution is unique and satisfies the conditions given later in Assumptions \ref{A1} and \ref{A2}.

Define $E_i \triangleq [0 \  E_{\delta_i} ]$ and  $R \triangleq [R_{\rom{r}} \ 0]$. Furthermore, let $e_i(t) \triangleq  y_i(t)-y_0(t)$ be the tracking error. We can then write the dynamics of each agent and its tracking error as
\begin{eqnarray}
\dot x_i(t) &=& A_ix_i(t)+B_iu_i(t) + E_i \omega(t), \quad x_i(0)=x_{i0}, \quad t \geq 0, \label{local_state_2}\\
e_i(t)&=& C_ix_i(t) + D_iu_i(t) - R\omega(t). \label{local_tracking_error_1}
\end{eqnarray} 
In this paper, the tracking error $e_i(t)$ is available to a nonempty proper subset of agents\footnote{\hspace{0.16cm}If all agents observe the leader, decentralized controllers can be designed for each agent even though the distributed controllers proposed here are still applicable.}. In particular, if node $v_i$ observes the leader node $v_0$, then there exists an edge $(v_0,v_i)$ with weighting gain $k_i>0$; otherwise $k_i = 0$. Each agent has also access to the relative output error; that is, $y_i(t)-y_j(t) $ for all $j \in N_i$. 
Similar to \cite{Farnaz_automatica}, the local virtual tracking error can be defined as 
\begin{eqnarray}
e_{\rom{v}i}(t)  &\triangleq&  \frac{1}{d_i+k_i} \Big[\hspace{-0.07cm}\sum_{j\in N_i}  a_{ij}\big(y_i(t) - y_j(t)\big) + k_i\big(y_i(t) - y_0(t)\big) \Big]. \label{virtual_local_error_1}
\end{eqnarray}

Now, we define three classes of distributed control laws based on additional available information to each agent:

\textbf{$1)$ Dynamic State Feedback.} If each agent has full access to its own state $x_i(t)$, then the dynamic state feedback control law is given by  
\begin{eqnarray}
u_i(t) &=& K_{1i}x_i(t) + K_{2i}z_i(t), \label{output_equation_of_dynamic_state_feedback_control} \\
\dot z_i(t) &=& G_{1i}z_i(t) + G_{2i}e_{\rom{v}i}(t),  \quad z_i(0)=z_{i0}, \quad t \geq 0, \label{state_equation_of_dynamic_state_feedback_control}
\end{eqnarray}
where $z_i(t) \in \IR^{n_{z_{1i}}}$ is the controller state and the quadruple $(K_{1i},K_{2i},G_{1i},G_{2i})$ is specified in Section \ref{first_controller}.

\textbf{$2)$ Dynamic Output Feedback with Local Measurement.} If each agent has local measurement output $y_{\rom{m}i}(t) \in\IR^{p_i} $ of the form
\begin{eqnarray}
y_{\rom{m}i}(t) &=& C_{\rom{m}i}x_i(t)+D_{\rom{m}i}u_i(t), \label{measured_output_1}
\end{eqnarray}	
then the dynamic output feedback control law with local measurement is given by  
\vspace{-0.3cm} 
\begin{eqnarray}
u_i(t) &=& \bar K_{i}z_i(t), \label{output_equation_of_dynamic_output_feedback_control_with_local_measurement} \\
\dot z_i(t) &=& M_{1i}z_i(t) + M_{2i}e_{\rom{v}i}(t)+ M_{3i}y_{\rom{m}i}(t),  \quad  z_i(0)=z_{i0}, \quad t \geq 0, \label{state_equation_of_dynamic_output_feedback_control_with_local_measurement}
\end{eqnarray}
where $z_i(t) \in \IR^{n_{z_{2i}}}$ is the controller state and the quadruple $(\bar K_{i},M_{1i},M_{2i},M_{3i})$ is specified in Section \ref{second_controller}. 

\textbf{$3)$ Dynamic Output Feedback.}	If each agent does not have additional information; that is, the local virtual tracking error $e_{\rom{v}i}(t)$ is the only available information to it, then the dynamic output feedback control law is given by 
\begin{eqnarray}
u_i(t) &=&  \bar K_{i}z_i(t), \label{output_equation_of_dynamic_output_feedback_control} \\
\dot z_i(t) &=&  M_{1i}z_i(t) +  M_{2i}e_{\rom{v}i}(t), \quad z_i(0)=z_{i0}, \quad t \geq 0, \label{state_equation_of_dynamic_output_feedback_control}
\end{eqnarray}
where $z_i(t) \in \IR^{n_{z_{2i}}}$ is the controller state and the triple $(\bar K_{i},M_{1i},M_{2i})$ is specified in Section \ref{third_controller}.

We now introduce the first and the second assumptions before defining the problem.

\assumption\label{A1}  $A_0 \in \IR^{q \times q}$ has no eigenvalues with negative real parts.
\vspace{-0.05cm}
\assumption\label{A2} There exists $\kappa > 0$ such that
\begin{eqnarray}
\norm{A_0 \omega(t) - \dot \omega(t) }_{2}  \leq \kappa < \infty, \quad \forall t \geq 0, \label{time_derivative_of_leader} \nonumber 
\end{eqnarray}
where $\dot \omega(t)$ is a piecewise continuous function\footnote{\hspace{0.16cm}We follow the definition given in page 650 of \cite{Khalil_book}.} of $t$.

Assumption \ref{A1} is standard in linear output regulation theory (e.g., see Remark 1.3 in \cite{Huang_book}). Assumption \ref{A2} is required to show the ultimate boundedness of the tracking error and
it automatically holds if the exosystem has an unforced linear time-invariant dynamics with the system matrix $A_0$. Note that these assumptions do not imply the exact knowledge of the exosystem. We refer to Remarks \ref{Remark2} and \ref{Remark3} for further discussions and Section \ref{INE} for illustrative examples on this point.

Based on the definition of the linear cooperative output regulation problem in  \cite{Su_2012,Huang_C_tac}, the problem considered in this paper is defined as follows.

\definition\label{Def2} Given the system in \eqref{local_state_2} and \eqref{local_tracking_error_1} together with the exosystem, which satisfies Assumptions \ref{A1} and \ref{A2}, and the fixed augmented directed graph $\mathcal{\bar G}$, find a distributed control law of the form \eqref{output_equation_of_dynamic_state_feedback_control} and \eqref{state_equation_of_dynamic_state_feedback_control}, or \eqref{output_equation_of_dynamic_output_feedback_control_with_local_measurement} and \eqref{state_equation_of_dynamic_output_feedback_control_with_local_measurement}, or 
\eqref{output_equation_of_dynamic_output_feedback_control} and \eqref{state_equation_of_dynamic_output_feedback_control} such that:

\indent $a)$ The resulting closed-loop system matrix  is Hurwitz. \\
\indent $b)$ The tracking error $e_i(t)$ is ultimately bounded with ultimate bound $b$ for all initial conditions of the closed-loop system and for all $i \in \mathcal{N}$; that is, there exists $b>0$ and for each initial condition of the closed-loop system, there is $T\geq 0$ such that $\norm{e_i(t)}_2 \leq  b, \ \forall t \geq T,  \ \forall i\in\mathcal{N}$. \\
\indent $c)$  If $\lim_{t\to\infty}  A_0 \omega(t) - \dot \omega(t)= 0$, then for all initial conditions of the closed-loop system \\  $\lim_{t\to\infty} e_i(t) = 0, \ \forall i\in\mathcal{N}$.

This paper makes the following additional assumptions to solve this problem.

\assumption\label{A3} The fixed augmented directed graph $\bar{\mathcal{G}}$ has a spanning tree with the root node being the leader node.
\vspace{-0.05cm}
\assumption\label{A4} The pair $(A_i,B_i)$ is stabilizable for all $i\in\mathcal{N}$.  
\vspace{-0.05cm}
\assumption\label{A5}
For all $\lambda \in \sigma(A_0)$, 
\begin{eqnarray} 
\rom{rank} 
\begin{bmatrix}
A_i-\lambda I_{n_i}  & B_i \\
C_i  & D_i  \\
\end{bmatrix} = n_i +p, \quad \forall i \in \mathcal{N}. \nonumber 
\end{eqnarray}
\vspace{-0.05cm}
\assumption\label{A6} As in \eqref{M1_M2_M3_Definition1_2}, the triple $(G_{1i},G_{2i}, 0)$ incorporates a $p$-copy internal model of $A_0$ for all $i\in\mathcal{N}$.
\vspace{-0.05cm}
\assumption\label{A7} The pair $(A_i,C_{\rom{m}i})$ is detectable for all $i\in\mathcal{N}$.
\vspace{-0.05cm}
\assumption\label{A8} The pair $(A_i,C_{i})$ is detectable for all $i\in\mathcal{N}$.  

Assumption \ref{A3} is natural to solve the stated problem (e.g., see Remark 3.2 in \cite{Lewis_book}). Similar to Assumption \ref{A1}, Assumptions \ref{A4}-\ref{A8} are standard in linear output regulation theory (e.g., see Chapter 1 of \cite{Huang_book}). We use Assumptions \ref{A1}-\ref{A6} for dynamic state feedback. To utilize some results from dynamic state feedback in the absence of full state information, each agent requires the estimation of its own state. For this purpose, Assumption \ref{A7} and Assumption \ref{A8} are included for dynamic output feedback with local measurement and dynamic output feedback, respectively.

\section{Solvability of the Problem}\label{DSF} 

For the three different distributed control laws introduced in Section \ref{PS}, this section investigates the solvability of the problem given in Definition \ref{Def2}. Specifically, the approach in this section is twofold. First, the property $a)$ of Definition \ref{Def2} is assumed and it is shown, under mild conditions, that the properties $b)$ and $c)$ of Definition \ref{Def2} are satisfied. Second, an agent-wise local sufficient condition (i.e., distributed criterion) is provided for the property $a)$ of Definition \ref{Def2} (i.e., the stability of the closed-loop system matrix) under standard assumptions.

Before studying the solvability of the problem for each distributed control law, we now present some definitions that are used throughout this section to express the closed-loop systems in compact forms, some results related to the communication graph topology, and a key lemma about the solvability of matrix equations, which play a crucial role on the solvability of the problem.

Define the following matrices:\\
 $ \Phi \triangleq \rom{diag} (\Phi_{1}, \ldots, \Phi _{N} ), \ \Phi = A,B,C,D,E$; 
$\Phi_{\rom{m}} \triangleq \rom{diag} (\Phi_{\rom{m}1}, \ldots, \Phi_{\rom{m}N} ), \ \Phi = C, D$;  \\
$ K_l  \triangleq \rom{diag} (K_{l1}, \ldots, K_{lN} ), \ l= 1,2 $; $A_{0\rom{a}} \triangleq I_N \otimes A_0 $, and $R_\rom{a} \triangleq I_N \otimes R $. \\
Further, let 
 $x(t) \triangleq [x_1^\rom{T}(t), \ldots, x_N^\rom{T}(t)]\mT$ $\in\IR^{\bar{n}}$, where $\bar{n} = \sum_{i=1}^{N} n_i $; $e(t) \triangleq [e_1^\rom{T}(t), \ldots, e_N^\rom{T}(t)]\mT$ $\in\IR^{Np}$, $e_{\rom{v}}(t) \triangleq [e_{\rom{v}1}^\rom{T}(t), \ldots, e_{\rom{v}N}^\rom{T}(t)]\mT$ $\in\IR^{Np}$, and $ {\omega_{{\rom{a}}}}(t) \triangleq \ \textbf{1}_N \otimes \omega(t) \in\IR^{Nq} $.

Observing $y_i(t)-y_j(t) = e_i(t)-e_j(t)$ and recalling $d_i = \sum_{j\in N_i} a_{ij} $, \eqref{virtual_local_error_1} can be equivalently written as
\begin{eqnarray}
e_{\rom{v}i}(t)&=& e_i(t) - \frac{1}{d_i+k_i} \sum_{j\in N_i} a_{ij}e_j(t). \label{virtual_local_error_2} 
\end{eqnarray}
Let $ \mathcal{F} \triangleq \rom{diag}\Big(\frac{1}{d_1+k_1}, \ldots, \frac{1}{d_N+k_N}\Big)$ and  $ \mathcal{W}  \triangleq (I_N-\mathcal{F}\mathcal{A}) \otimes I_p$. Here, it should be noted that $d_i+k_i >0,\ \forall i \in \mathcal{N}$ by Assumption \ref{A3};  hence, $\mathcal{F}$ is well-defined. From \eqref{virtual_local_error_2}, we have 
\begin{eqnarray}
e_{\rom{v}}(t)&=& \mathcal{W}e(t). \label{e_v}  
\end{eqnarray}

Similar to Lemma 3.3 in \cite{Lewis_book}, we next present the following lemma for $I_N-\mathcal{F}\mathcal{A}$.

\lemma\label{Lemma1} Under Assumption \ref{A3}, $I_N-\mathcal{F}\mathcal{A}$ is nonsingular. In addition, all its eigenvalues have positive real parts.

\textit{Proof.} Under Assumption \ref{A3}, $I_N-\mathcal{F}\mathcal{A}$ satisfies the conditions of the theorem in \cite{Shivakumar}. Thus,  it is nonsingular.
Since the singularity is eliminated, all the eigenvalues of $I_N-\mathcal{F}\mathcal{A}$ have positive real parts by the Gershgorin circle theorem  (e.g., see Fact 4.10.17 in \cite{Bernstein_book}). \hfill$\blacksquare$ 

\remark\label{Remark1} Since $I_N-\mathcal{F}\mathcal{A}$ is nonsingular under Assumption \ref{A3}, so is   $\mathcal{W}$ by Proposition 7.1.7 in \cite{Bernstein_book}. Then, it is clear from \eqref{e_v} that 
$e_i(t)$ is bounded for all $i \in \mathcal{N}$ if and only if $e_{\rom{v}i}(t)$ is bounded for all $i \in \mathcal{N}$;  $\lim_{t\to\infty} e_i(t) = 0, \ \forall i\in\mathcal{N}$  if and only if $\lim_{t\to\infty} e_{\rom{v}i}(t) = 0, \ \forall i\in\mathcal{N}$.

We now investigate the spectral radius of $\mathcal{FA}$.

\lemma\label{Lemma2} Under Assumption \ref{A3}, $\rho(\mathcal{FA})<1$.

\textit{Proof.} By Lemma \ref{Lemma1}, all the eigenvalues of $I_N-\mathcal{F}\mathcal{A}$  have positive real parts under Assumption \ref{A3}. This directly implies from Fact 6.2.1.4 in \cite{Vidyasagar_1981} that the leading principal minors of $I_N-\mathcal{F}\mathcal{A}$ are all positive as $I_N-\mathcal{F}\mathcal{A}$ is a square matrix whose off-diagonal elements are all nonpositive. Since $\mathcal{FA}$ is a nonnegative square matrix and the leading principal minors of $I_N-\mathcal{F}\mathcal{A}$ are all positive, $\rho(\mathcal{FA})<1$ from Lemma 6.2.1.8 in \cite{Vidyasagar_1981}. \hfill$\blacksquare$  

Finally, we introduce the key lemma that extends the field of application of Lemma 1.27 in \cite{Huang_book} to heterogeneous (in dynamics and dimension) linear time-invariant multiagent systems over general fixed directed communication graph topologies.

\lemma\label{Lemma3}\hspace{-0.29cm}\footnote{\hspace{0.16cm}To investigate the solvability of a matrix equation that is obtained for a different problem setting with the distributed dynamic state feedback control law, the authors of \cite{Wang_Chen_Liu_2018} utilized the same logic in the proof of Lemma \ref{Lemma3} (see Section 3.1 in \cite{Wang_Chen_Liu_2018}). }
Let Assumptions \ref{A1} and \ref{A3} hold. Suppose the triple $(M_1,M_2,M_3)$ incorporates an $Np$-copy internal model of $A_{0\rom{a}}$.
If 
\begin{eqnarray}
A_{\rom{c}}  &\triangleq&  \begin{bmatrix}
\hat A & \hat B \\
M_2\mathcal{W}\hat C+M_3 \hat C_{\rom{m}} \ & \ M_1 +M_2\mathcal{W}\hat D + M_3 \hat D_{\rom{m}} 
\end{bmatrix} \nonumber 
\end{eqnarray}	
is Hurwitz, where $\hat A$, $\hat B$, $\hat C$, $\hat C_{\rom{m}}$, $\hat D$, and $ \hat D_{\rom{m}}$ are any matrices with appropriate dimensions, then the matrix equations
\begin{eqnarray}
XA_{0\rom{a}}  &=& \hat A X + \hat B Z + \hat E, \label{lme1_generic_state} \\
ZA_{0\rom{a}}  &=& M_1 Z + M_2 \mathcal{W}(\hat C X + \hat DZ + \hat F) + M_3(\hat C_{\rom{m}} X + \hat D_{\rom{m}}Z) \label{lme2_generic_auxiliary},
\end{eqnarray}
have unique solutions $X$ and $Z$ for any matrices $\hat E$ and $\hat F$ of appropriate dimensions. Furthermore,  $X$ and $Z$ satisfy
\begin{eqnarray}
0 &=& \hat C X + \hat DZ +\hat F \label{lme3_generic_error}.
\end{eqnarray}
In other words, the conclusion is that the matrix equations
\begin{eqnarray}
X_{\rom{c}}A_{0\rom{a}} &=& A_{\rom{c}}X_{\rom{c}} + B_{\rom{c}}, \label{lme1_generic_state_and_auxiliary} \\
0 &=& C_{\rom{c}}X_{\rom{c}} + D_{\rom{c}}, \label{lme2_generic_error}
\end{eqnarray}
have a unique solution $X_{\rom{c}}$, where 
\begin{eqnarray}
X_{\rom{c}}  &=& \begin{bmatrix}
X \\
Z 
\end{bmatrix}, \ \
B_{\rom{c}}  =  \begin{bmatrix}
\hat E  \\
M_2 \mathcal{W} \hat F
\end{bmatrix}, \ 
C_{\rom{c}}   =  \begin{bmatrix}
\hat C & \hat D  
\end{bmatrix}, \ D_{\rom{c}}   =  \hat F. \nonumber 
\end{eqnarray}	

\textit{Proof.} Note that \eqref{lme1_generic_state} and \eqref{lme2_generic_auxiliary} (respectively, \eqref{lme3_generic_error}) can be equivalently written as \eqref{lme1_generic_state_and_auxiliary} (respectively, \eqref{lme2_generic_error}). Note also that  $\sigma (A_{0\rom{a}}) = \sigma(A_0)$. Since Assumption \ref{A1} holds and $A_{\rom{c}}$ is Hurwitz, $A_{0\rom{a}}$ and $A_{\rom{c}}$ have no eigenvalues in common. Thus, the Sylvester equation in \eqref{lme1_generic_state_and_auxiliary} has a unique solution $X_{\rom{c}} = [X^{\rom{T}} \ Z^\rom{T}]\mT $ by the first part of Proposition A.2 in \cite{Huang_book}. 
In addition, we show that $X$ and $Z$ also satisfy \eqref{lme3_generic_error}. To this end, let $ \bar \gamma \triangleq \hat C X + \hat DZ +\hat F  $.  Since the triple $(M_1,M_2,M_3)$ incorporates an $Np$-copy internal model of $A_{0\rom{a}}$, it has the form given by \eqref{M1_M2_M3_Definition1_1} or \eqref{M1_M2_M3_Definition1_2}. 
If it takes the form \eqref{M1_M2_M3_Definition1_1}, let $[{\hat \theta}^{\rom{T}} \ {\bar \theta}^\rom{T}]\mT \triangleq T^{-1}Z$, where $\bar \theta$ has as many rows as those of $G_1$.
Premultiplying \eqref{lme2_generic_auxiliary} by $T^{-1}$ and using the
foregoing definitions, we obtain 
\begin{eqnarray}
\bar \theta  A_{0\rom{a}} &=& G_1 \bar \theta + G_2 \mathcal{W} \bar \gamma. \label{prekey_equation}
\end{eqnarray}
Note that if the triple $(M_1,M_2,M_3)$ takes the form \eqref{M1_M2_M3_Definition1_2}, \eqref{lme2_generic_auxiliary} already satisfies \eqref{prekey_equation}, where $\bar \theta = Z$.
Let  $\gamma \triangleq \mathcal{W} \bar \gamma$; then, \eqref{prekey_equation} is in the form of (1.74) in \cite{Huang_book}. Hence, $\gamma = 0$ by the proof of Lemma 1.27 in \cite{Huang_book}. We know from Remark \ref{Remark1} that  $\mathcal{W}$ is nonsingular under Assumption \ref{A3}. As a consequence, $\gamma = 0$ implies $\bar \gamma=0$. This completes the proof of this lemma. 
\hfill$\blacksquare$   

\subsection{Dynamic State Feedback}\label{first_controller}

Let $z(t) \triangleq [z_1^\rom{T}(t), \ldots, z_N^\rom{T}(t)]\mT   \in \IR^{\bar n_{z_{1}}}$, where $\bar n_{z_{1}}= \sum_{i=1}^{N} n_{z_{1i}}$, and $ G_l  \triangleq \rom{diag} (G_{l1}, \ldots, G_{lN} ), \ l= 1,2 $. Inserting  \eqref{output_equation_of_dynamic_state_feedback_control} into \eqref{local_state_2} and \eqref{local_tracking_error_1}, and using the above definitions, \eqref{local_state_2}, \eqref{state_equation_of_dynamic_state_feedback_control}, and \eqref{local_tracking_error_1} can be compactly written as 
\begin{eqnarray}
\dot x(t) &=& (A+BK_1)x(t)+ BK_2z(t) + E{\omega_{{\rom{a}}}}(t) ,  \quad x(0)=x_{0},  \quad t \geq 0, \label{global_state_1}\\
\dot z(t) &=& G_1z(t)+G_2 e_{\rom{v}}(t), \quad \hspace{2.96cm} z(0)=z_0,  \quad t \geq 0, \label{global_controller_state_1} \\
e(t)&=& (C+DK_1)x(t) + DK_2z(t) - R_{\rom{a}}{\omega_{{\rom{a}}}}(t). \label{global_tracking_error_1}
\end{eqnarray} 
Next, insert \eqref{global_tracking_error_1} into \eqref{e_v} and replace the obtained expression with the one in \eqref{global_controller_state_1}. Define $x_{\rom{g}}(t)  \triangleq [x^\rom{T}(t), z^\rom{T}(t)]\mT \in \IR^{\bar n +\bar n_{z_{1}}}$. Then, the closed-loop system of \eqref{local_state_2}-\eqref{state_equation_of_dynamic_state_feedback_control} becomes
\begin{eqnarray}
\dot x_{\rom{g}}(t) &=& A_\rom{g} x_{\rom{g}}(t)+B_\rom{g} {\omega_{{\rom{a}}}}(t),  \quad x_{\rom{g}}(0)=x_{\rom{g}0}, \quad t  \geq 0, \label{full_state_global_state}\\
e(t)&=& C_\rom{g}x_{\rom{g}}(t)+ D_\rom{g} {\omega_{{\rom{a}}}}(t), \label{full_state_global_error}
\end{eqnarray}
where 
\begin{eqnarray}
A_{\rom{g}}  &=&  \begin{bmatrix}
A+BK_1  & BK_2 \\
G_2\mathcal{W}(C+DK_1)  & G_1 +G_2\mathcal{W}DK_2 
\end{bmatrix}, \ 
B_{\rom{g}}  =  \begin{bmatrix}
E    \\
-G_2\mathcal{W} R_{\rom{a}} 
\end{bmatrix},   \nonumber  \\
C_{\rom{g}} &=&  \begin{bmatrix}
C+DK_1 & DK_2  
\end{bmatrix}, \ D_{\rom{g}} = - R_{\rom{a}}. \nonumber 
\end{eqnarray}	

\theorem\label{Theorem1} Let Assumptions \ref{A1}-\ref{A3} and \ref{A6} hold. If $A_{\rom{g}}$ is Hurwitz, then the distributed dynamic state feedback control given by \eqref{output_equation_of_dynamic_state_feedback_control} and \eqref{state_equation_of_dynamic_state_feedback_control} solves the problem in Definition \ref{Def2}.

\textit{Proof.} By the definition of $A_{0\rom{a}}$, the minimal polynomials for $A_{0\rom{a}}$ and $A_0$ are the same. Thus, the triple $(G_1,G_2,0)$ incorporates an $Np$-copy internal model of $A_{0\rom{a}}$ under Assumption \ref{A6}. Let $(M_1,M_2,M_3) \triangleq (G_1,G_2,0)$.
Let also $\hat A \triangleq A+BK_1$, $\hat B \triangleq BK_2$, $\hat C \triangleq C+DK_1$, $\hat C_{\rom{m}} \triangleq 0 $, $\hat D \triangleq DK_2$, $\hat D_{\rom{m}} \triangleq 0 $, $\hat E \triangleq E$, and $\hat F \triangleq - R_{\rom{a}} $. Then, the quadruple $(A_{\rom{g}},B_{\rom{g}},C_{\rom{g}},D_{\rom{g}})$ takes the form of $(A_{\rom{c}},B_{\rom{c}},C_{\rom{c}},D_{\rom{c}})$ in Lemma \ref{Lemma3}. In addition, $A_{\rom{g}}$ is Hurwitz and  Assumptions \ref{A1} and \ref{A3} hold. Hence, Lemma \ref{Lemma3} is applicable and it implies that the matrix equations
\begin{eqnarray}
X_{\rom{g}}A_{0\rom{a}} &=& A_{\rom{g}}X_{\rom{g}} + B_{\rom{g}}, \label{lme1_state_feedback} \\
0 &=& C_{\rom{g}}X_{\rom{g}} + D_{\rom{g}}, \label{lme2_state_feedback}
\end{eqnarray}
have a unique solution $X_{\rom{g}}$. We also refer to Appendix A for additional discussions on the solvability of \eqref{lme1_state_feedback} and \eqref{lme2_state_feedback}.

Under Assumption \ref{A2}, $\norm{A_{0\rom{a}}  {\omega_{{\rom{a}}}}(t) - {\dot \omega_{{\rom{a}}}}(t)}_2 \leq \sqrt{N}\kappa, \ \forall t  \geq 0 $ since $\norm{A_{0\rom{a}}  {\omega_{{\rom{a}}}}(t) - {\dot \omega_{{\rom{a}}}}(t)}_2^2 = N \norm{A_0 \omega(t) \\- \dot \omega(t) }_{2}^2 $. Let $\bar x_{\rom{g}} (t) \triangleq x_{\rom{g}}(t) - X_{\rom{g}}{\omega_{{\rom{a}}}}(t) $. Then,  using the definition of $ \bar x_{\rom{g}} (t)$  and \eqref{lme1_state_feedback} and \eqref{lme2_state_feedback}, we can rewrite \eqref{full_state_global_state} and \eqref{full_state_global_error} as 
\begin{eqnarray}
\dot {\bar x}_{\rom{g}}(t) &=& A_\rom{g} {\bar x}_{\rom{g}}(t) + X_{\rom{g}}(A_{0\rom{a}}  {\omega_{{\rom{a}}}}(t) - {\dot \omega_{{\rom{a}}}}(t) ), \quad \bar x_{\rom{g}}(0)= \bar x_{\rom{g}0}, \quad t  \geq 0, \label{full_state_global_assistant_state}  \\
e(t)&=& C_\rom{g} \bar x_{\rom{g}}(t). \label{full_state_global_error_due_to_assitant}
\end{eqnarray}
Now, the solution of \eqref{full_state_global_assistant_state} can be written as 
\begin{eqnarray}
{\bar x}_{\rom{g}}(t) =  e^{A_\rom{g}t}\bar x_{\rom{g}0} + \int_{0}^{t}  e^{A_\rom{g}(t-\tau)} X_{\rom{g}}(A_{0\rom{a}}  {\omega_{{\rom{a}}}}(\tau) - {\dot \omega_{{\rom{a}}}}(\tau) ) \rom{d}\tau. \label{full_state_assistant_solution} \nonumber 
\end{eqnarray}
Since $A_{\rom{g}}$ is Hurwitz, there exist $c>0$ and $\alpha>0$ such that $\norm{e^{A_\rom{g}t}}_2 \leq ce^{-\alpha t}, \ \forall t \geq 0$ (e.g., see Lecture 8.3 in \cite{Hespanha_book}). Owing to this bound and the bound on $\norm{A_{0\rom{a}}  {\omega_{{\rom{a}}}}(t) - {\dot \omega_{{\rom{a}}}}(t)}_2$, we have the following inequality
\begin{eqnarray}
\norm{{\bar x}_{\rom{g}}(t) }_2  \leq  ce^{-\alpha t}\norm{\bar x_{\rom{g}0}}_2 +  \frac{c\norm{X_\rom{g}}_2 }{\alpha}  \sqrt{N}\kappa, \quad \forall t  \geq 0.  \label{full_state_assistant_inequality}  \nonumber
\end{eqnarray}
Using the fact $\norm{e_i(t)}_2 \leq \norm{e(t)}_2, \  \forall  i\in\mathcal{N}$ and observing $\norm{e(t)}_2 \leq \norm{C_\rom{g}}_2 \norm{\bar x_{\rom{g}}(t)}_2$ from \eqref{full_state_global_error_due_to_assitant}, we arrive
\begin{eqnarray}
\norm{e_i(t)}_2  \leq   ce^{-\alpha t}\norm{C_\rom{g}}_2 \norm{\bar x_{\rom{g}0}}_2 +  b' ,  \quad \forall t  \geq 0, \quad \forall i \in \mathcal{N},  \label{tracking_error_inequality} \nonumber 
\end{eqnarray}
where $b' = c\norm{C_\rom{g}}_2 \norm{X_\rom{g}}_2\sqrt{N}\kappa  \alpha^{-1} $. For a given $\epsilon>0$, we have either  $c\norm{C_\rom{g}}_2 \norm{\bar x_{\rom{g}0}}_2> \epsilon$ or \\ $c\norm{C_\rom{g}}_2 \norm{\bar x_{\rom{g}0}}_2 \leq \epsilon$. In the former case, it can be readily shown that $ce^{-\alpha t}\norm{C_\rom{g}}_2 \norm{\bar x_{\rom{g}0}}_2 \leq \epsilon, \ \forall t \geq T$ with $T = \alpha^{-1} \rom {ln}\Big(\frac{c\norm{C_\rom{g}}_2 \norm{\bar x_{\rom{g}0}}_2}{\epsilon}\Big)  > 0$. In the latter case, the foregoing inequality trivially holds for all $t \geq 0$. Thus,  $e_i(t)$ is ultimately bounded with the ultimate bound $b \triangleq b'+\epsilon$ for all $\bar x_{\rom{g}0} $, which is also true for all $ x_{\rom{g}0} $, and 
for all $i \in \mathcal{N}$.

If $\lim_{t\to\infty}  A_0 \omega(t) - \dot \omega(t)= 0$, then $ \lim_{t\to\infty} A_{0\rom{a}}  {\omega_{{\rom{a}}}}(t) - {\dot \omega_{{\rom{a}}}}(t) = 0 $. Since  $A_{\rom{g}}$ is Hurwitz  and the system in \eqref{full_state_global_assistant_state} is linear time-invariant when $A_{0\rom{a}}  {\omega_{{\rom{a}}}}(t) - {\dot \omega_{{\rom{a}}}}(t)$  is viewed as an input to the system,  \eqref{full_state_global_assistant_state} is input-to-state stable with respect to this piecewise continuous input (e.g., see Chapter 4.9 in \cite{Khalil_book}).
Thus, $ \lim_{t\to\infty} A_{0\rom{a}}  {\omega_{{\rom{a}}}}(t) - {\dot \omega_{{\rom{a}}}}(t) = 0 $ implies $\lim_{t\to\infty} \bar x_{\rom{g}}(t) = 0 $ for all $\bar x_{\rom{g}0} $ (e.g., see Exercise 4.58 in \cite{Khalil_book}). Finally, it follows from \eqref{full_state_global_error_due_to_assitant} that for all $ x_{\rom{g}0} $ $\lim_{t\to\infty} e_i(t)= 0, \ \forall i \in \mathcal{N} $. \hfill$\blacksquare$ 

\remark\label{Remark2} The ultimate bound $b$ of the tracking error for each agent is associated with the bound $\kappa$ in Assumption \ref{A2}. Specifically, as $\kappa$ decreases (respectively, increases), $b$ decreases (respectively, increases). To elucidate the role of Assumptions \ref{A1} and \ref{A2} in practice, we consider the following possible scenarios:

 $a)$
When the piecewise continuity and boundedness of $\dot \omega(t)$ are the only information that is available to a control designer, the triple $(0,I_p,0)$ incorporating a $p$-copy internal model of $A_0 = 0$ is quite natural; hence, \eqref{state_equation_of_dynamic_state_feedback_control} becomes a distributed integrator. Moreover, $X_{\rom{g}}$ in $b$ can be explicitly expressed in terms of $A_{\rom{g}}$ and  $B_{\rom{g}}$; that is,  $X_{\rom{g}} =- A_{\rom{g}}^{-1} B_{\rom{g}}$  by \eqref{lme1_state_feedback}.

$b)$ When the piecewise continuity and boundedness of $\dot \omega(t)$, the boundedness of $ \omega(t)$, and some frequencies in $ \omega(t)$ are available to a control designer, the triple $(G_{1i},G_{2i},0)$ incorporating a $p$-copy internal model of $A_0$, which includes these frequencies and zero eigenvalues, is an alternative to the pure distributed integrator. 

\remark\label{Remark3} As it is shown in Theorem \ref{Theorem1}, asymptotic synchronization is achieved when \\ $\lim_{t\to\infty}  A_0 \omega(t) - \dot \omega(t)= 0$. We now provide sufficient conditions to check this condition as follows\footnote{\hspace{0.16cm}If $A_0 = 0$, one should read $\lim_{t\to\infty}  \dot \omega(t)= 0$ in place of $\lim_{t\to\infty}  A_0 \omega(t) - \dot \omega(t)= 0$; hence,  $\omega(t) \equiv \omega^\star$  ($\omega^\star$ is finite) in place of $a)$, and $\lim_{t\to\infty}   \omega(t)= \omega^\star$ and $\dot \omega(t)$ is uniformly continuous on $[0, \infty)$ in place of $b)$.}. If one of the following conditions holds

$a)$ $\dot \omega(t) = A_0 \omega(t), \quad \omega(0)=\omega_0, \quad t\geq 0$;

$b)$ $\lim_{t\to\infty} e^{A_0t}\omega_0-\omega(t) = 0 $, where $\omega_0 = \omega(0)$, and $A_0 e^{A_0t}\omega_0 -\dot \omega(t) $ is uniformly continuous on $[0, \infty)$,

\noindent then $\lim_{t\to\infty}  A_0 \omega(t) - \dot \omega(t)= 0$. Note that $a)$ clearly implies $b)$. From Barbalat's lemma given by Lemma 8.2 in \cite{Lavretsky_Wise_book}, $b)$ implies that  $\lim_{t\to\infty} A_0 e^{A_0t}\omega_0 -\dot \omega(t) = 0 $. Thus, $\lim_{t\to\infty}  A_0 \omega(t) - \dot \omega(t)= A_0 \lim_{t\to\infty} \omega(t) - e^{A_0t}\omega_0  + \lim_{t\to\infty} A_0 e^{A_0t}\omega_0 -\dot \omega(t) = 0 $. In general, asymptotic synchronization results in the literature (e.g., see \cite{Su_2012,Huang_C_tac,Farnaz_automatica}) are obtained under the condition $a)$. It is clear that this paper covers all class of functions generated under the condition $a)$. 

To obtain an agent-wise local sufficient condition assuring the property $a)$ of Definition \ref{Def2} under some standard assumptions, let $\xi_i (t)\triangleq [x_i^\rom{T}(t), z_i^\rom{T}(t)]^\rom{T} \in \IR^{n_i + n_{z_{1i}}} $, $\mu_i(t) \triangleq   \frac{1}{d_i + k_i}\sum_{j\in N_i} a_{ij}e_j(t)$,
\begin{eqnarray}
\bar A_i &\triangleq &  \begin{bmatrix}
A_i & 0 \\
G_{2i}C_i  & G_{1i} 
\end{bmatrix}, \ 
\bar B_{i} \triangleq   \begin{bmatrix}
B_i  \\
G_{2i}D_i  
\end{bmatrix}, \  B_{\rom{f}i} \triangleq
\begin{bmatrix}
0   \\
-G_{2i} 
\end{bmatrix}, \nonumber 
\end{eqnarray}	
and $\bar C_i \triangleq [C_i \ 0]$. Furthermore, consider \eqref{local_state_2}, \eqref{state_equation_of_dynamic_state_feedback_control}, \eqref{virtual_local_error_2}, and \eqref{local_tracking_error_1} when $\omega(t)  \equiv 0$. We now have
\begin{eqnarray}
\dot \xi_i(t) &=&\bar A_i \xi_i(t)+ \bar B_i u_i(t) + B_{\rom{f}i} \mu_i(t),   \quad \xi_i(0)=\xi_{i0},  \quad t  \geq 0, \label{local_open_loop_state} \\
e_i(t) &=& \bar C_i \xi_i(t) +D_i u_i(t). \label{local_open_loop_output}
\end{eqnarray}
Next, define the matrices
\begin{eqnarray}
A_{\rom{f}i} &\triangleq&  \begin{bmatrix}
A_i+B_iK_{1i}  & B_iK_{2i} \\
G_{2i}(C_i+D_iK_{1i})  & G_{1i} +G_{2i}D_iK_{2i} 
\end{bmatrix}, \nonumber  \\  
C_{\rom{f}i} &\triangleq&  \begin{bmatrix}
C_i+D_iK_{1i} & D_iK_{2i}  
\end{bmatrix}. \nonumber
\end{eqnarray}
Using \eqref{output_equation_of_dynamic_state_feedback_control}, \eqref{local_open_loop_state} and \eqref{local_open_loop_output} can be written as
\begin{eqnarray}
\dot \xi_i(t) &=& A_{\rom{f}i} \xi_i(t) + B_{\rom{f}i}  \mu_i(t),  \quad \xi_i(0)=\xi_{i0}, \quad t \geq 0, \label{local_closed_loop_state} \\
e_i(t)&=& C_{\rom{f}i} \xi_i(t). \label{local_closed_loop_output}
\end{eqnarray}

Let, in addition, $ \Psi_{\rom{f}} \triangleq \rom{diag} (\Psi_{\rom{f}1}, \ldots, \Psi_{\rom{f}N})$,  $\Psi = A,B,C$ and $\xi(t)  \triangleq [\xi_1^\rom{T}(t), \ldots, \xi_N^\rom{T}(t)]\mT $. Then, \eqref{local_closed_loop_state}  and \eqref{local_closed_loop_output} can be put into the compact form given by 
\vspace{-0.3cm}
\begin{eqnarray}
\dot \xi(t) &=& A_{\rom{f}} \xi(t) + B_{\rom{f}}(\mathcal{F}\mathcal{A} \otimes I_p  ) \tilde w(t),  \quad \xi(0)=\xi_{0},  \quad t \geq 0, \label{global_closed_loop_state} \\
\tilde z(t)&=& C_{\rom{f}} \xi(t), \label{global_closed_loop_output}
\end{eqnarray}
where $e(t) = \tilde w(t) = \tilde z(t)$. Observe that the system in \eqref{global_closed_loop_state} and \eqref{global_closed_loop_output} takes the form of (12) in \cite{Huang_C_tac}. Therefore, one may think of resorting Theorem 2 in \cite{Huang_C_tac} at first sight. However, the statement of Theorem 2 in  \cite{Huang_C_tac} is not correct as it is written; we refer to Appendix B for a counterexample.

\textit{This paragraph uses the notation and the terminology from \cite{Huang_C_tac}. Readers are referred to (12), Theorem 1, Theorem  2, and Lemma 8 in \cite{Huang_C_tac}.} It should be noted that Theorem 2 relies on Theorem 1 and this theorem is derived by means of Theorem 11.8 and Lemma 11.2 in \cite{Zhou_Doyle_Glover_book}. 
According to the mentioned results and Chapter 5.3, which is devoted to the notion of internal stability for the system of interest, in \cite{Zhou_Doyle_Glover_book}, it is clear that the following condition should be added to the hypotheses of Theorem 1: \textit{Let the realization of $T(s)$ given by (12) be stabilizable and detectable.} With this modification, not only the theoretical gap in Theorem 1 but also  the one in Theorem 2 is filled. 
However, a simple point in the proof of Theorem 2 still needs to be clarified. The spectral radius of $\tilde T (j\omega)$ in the proof of Theorem 2 is upper bounded by applying Lemma 8.
Since Lemma 8 is applied, we infer that $\rom{diag}(\norm{T_1(j\omega)},  \ldots, \norm{T_N(j\omega)} )$ is regarded as a positive definite diagonal matrix, but its proof is not given. The foregoing diagonal matrix is necessarily positive semidefinite; hence, we only question\footnote{\hspace{0.16cm}Considering Kalman decomposition  (e.g., see Theorem 16.3 in \cite{Hespanha_book}), one can easily construct a linear time-invariant system with Hurwitz system matrix, nonzero input and output matrices, and zero direct feedfeedthrough matrix such that its transfer matrix is zero.}  whether $T_i(s)=0$ for some $i$. Instead of investigating the corresponding realizations, we extend Lemma 8 to positive semidefinite diagonal matrices as follows.

\lemma\label{Lemma4} Let $Q \in \IR^{n \times n}$ be a nonnegative matrix. If $\Lambda \in \IR^{n \times n}$ is a positive semidefinite diagonal matrix, then $\rho(\Lambda Q) \leq \rho(\Lambda) \rho(Q)$.

\textit{Proof.} Let $\Lambda \triangleq \rom{diag}(\lambda_1, \ldots, \lambda_n)$ be positive semidefinite. If $\Lambda = 0$, the inequality holds trivially. We therefore assume that there exists a $\lambda_i >0$ for some $i$; hence, $\rho(\Lambda)>0$. Let $\bar \Lambda \triangleq \rom{diag}(\bar \lambda_1, \ldots, \bar \lambda_n) $, where $\bar \lambda_i = \rho(\Lambda)$ if $\lambda_i =0$, $\bar \lambda_i = \lambda_i$ otherwise. By construction,  $\Lambda \leq \bar \Lambda$, $\rho(\Lambda)=\rho(\bar \Lambda)$, and $\bar \Lambda $ is a positive definite diagonal matrix.
Since  $\Lambda \leq \bar \Lambda $ and $Q$ is nonnegative, $\Lambda Q \leq \bar \Lambda Q$. By the corollary in page 27 of \cite{Nonnegative_1994}, $\rho( \Lambda Q) \leq \rho(\bar \Lambda Q)$. 
 Applying  Lemma 8 in \cite{Huang_C_tac} to $\bar \Lambda Q$, we also have  $\rho(\bar \Lambda Q) \leq \rho(\bar \Lambda) \rho(Q)$. Since $\rho(\Lambda)=\rho(\bar \Lambda)$, we establish the desired inequality.
\hfill$\blacksquare$
 
It is well known that the system in \eqref{global_closed_loop_state} and \eqref{global_closed_loop_output} is stabilizable and detectable if $A_{\rom{f}}$ is Hurwitz. Thus, the new condition is satisfied if $A_{\rom{f}i}$ is Hurwitz for all $i\in\mathcal{N}$. 

\remark\label{Remark4} Assumptions \ref{A4}-\ref{A6} ensure the stabilizability of the pair $(\bar A_i,\bar B_i)$ for all $i\in\mathcal{N}$ by Lemma 1.26 in \cite{Huang_book}. Therefore, $K_{1i}$ and $K_{2i}$ can always be chosen such that $A_{\rom{f}i}$ is Hurwitz for all $i\in\mathcal{N}$.

Let $g_{\rom{f}i}(s) \triangleq C_{\rom{f}i}(sI -A_{\rom{f}i})^{-1}B_{\rom{f}i}$. We now state the following theorem for the dynamic state feedback case.

\theorem\label{Theorem2} Let Assumption \ref{A3} hold and $A_{\rom{f}i}$ be Hurwitz for all $i \in\mathcal{N}$. If 
\begin{eqnarray}
\norm {g_{\rom{f}i}}_\infty \ \hspace{-0.05cm} \rho(\mathcal{F}\mathcal{A}) < 1, \quad \forall i \in \mathcal{N}, \label{local_sufficient_cond}
\end{eqnarray}
where $\norm {g_{\rom{f}i}}_\infty$ is the $H_\infty$ norm of $g_{\rom{f}i}(s)$, then $A_{\rom{g}}$ is Hurwitz.

\textit{Proof.} It follows from Theorem 2 in \cite{Huang_C_tac}  and the above discussion. \hfill$\blacksquare$ 

\remark\label{Remark5} The inequality given by \eqref{local_sufficient_cond} is an agent-wise local sufficient condition; that is,
it paves the way for independent controller design for each agent.
For the connection between this condition and an algebraic Riccati equation (respectively, linear matrix inequality), we refer to Lemma 9 in \cite{Huang_C_tac} (respectively, Theorem 6 in \cite{Farnaz_automatica}). 
Moreover, we know from Lemma \ref{Lemma2} that $\rho(\mathcal{F}\mathcal{A}) < 1$ under Assumption \ref{A3}. 
Therefore, we can restate Theorem \ref{Theorem2} by replacing \eqref{local_sufficient_cond} with $\norm {g_{\rom{f}i}}_\infty \leq 1, \ \forall i \in \mathcal{N}$. In this statement, although the condition becomes more conservative, it is not only agent-wise local but also graph-wise local except Assumption \ref{A3}. Finally, it should be noted that if the graph $\mathcal{G}$ considered in Theorem \ref{Theorem2} contains no loop (i.e., acyclic), then the nodes in $\mathcal{G}$ can be relabelled such that $i > j$ when $(v_j, v_i) \in \mathcal{E}$. Thus, $\mathcal{A}$ is similar to a lower triangular matrix with zero diagonal entries, so is $\mathcal{FA}$. This implies that  $\rho(\mathcal{FA}) = 0$; hence, Theorem \ref{Theorem2} does not require the condition given by \eqref{local_sufficient_cond} anymore. In terms of being agent-wise and graph-wise local, this special case is consistent with the result in \cite{Wang_tac_2010}.


\subsection{Dynamic Output Feedback with Local Measurement} \label{second_controller}

Let $z_i(t) \triangleq [\hat x_i^\rom{T}(t), \bar z_i^\rom{T}(t)]\mT   \in \IR^{n_{z_{2i}}}$, where $\hat x_i(t)$ is the estimate of the state $x_i(t)$, $\bar K_i \triangleq [K_{1i} \ K_{2i}]$, and \eqref{output_equation_of_dynamic_output_feedback_control_with_local_measurement} have the form given by
\begin{eqnarray}
u_i(t) &=& K_{1i}\hat x_i(t) + K_{2i} \bar z_i(t). \label{output_equation_of_dynamic_output_feedback_control_with_local_measurement_specified}
\end{eqnarray}
To estimate the state $x_i(t)$, the following local Luenberger observer is employed
\begin{eqnarray}
\dot {\hat x}_i(t)  &=&   A_i{\hat x}_i(t) +  B_iu_i(t)  +   H_i \big(y_{\rom{m}i}(t)-   C_{\rom{m}i}{\hat x}_i(t)  -   D_{\rom{m}i}u_i(t) \big), \quad {\hat x}_i(0)={\hat x}_{i0},  \quad t \geq 0, \label{local_state_observer_1}
\end{eqnarray} 
where $H_i$ is the observer gain matrix. 
Using \eqref{output_equation_of_dynamic_output_feedback_control_with_local_measurement_specified}, we can write \eqref{local_state_observer_1} as 
\begin{eqnarray}
\dot {\hat x}_i(t)  &=& \big(A_i+B_i K_{1i}-H_i(C_{\rom{m}i}+D_{\rom{m}i}K_{1i})\big){\hat x}_i(t) +H_iy_{\rom{m}i}(t) +  (B_i-H_iD_{\rom{m}i})K_{2i}\bar z_i(t),  \nonumber \\ && \hspace{8.9cm}  \quad {\hat x}_i(0)={\hat x}_{i0},  \quad t \geq 0.  \label{local_state_observer_2}
\end{eqnarray} 
Let also $\bar z_i(t)$ evolve according to the dynamics given by
\begin{eqnarray}
\dot {\bar z}_i(t) &=& G_{1i}\bar z_i(t) + G_{2i}e_{\rom{v}i}(t),  \quad \bar z_i(0)= \bar z_{i0}, \quad t \geq 0. \label{partial_state_equation_of_dynamic_output_feedback_control_with_local_measurement}
\end{eqnarray}
By \eqref{local_state_observer_2} and \eqref{partial_state_equation_of_dynamic_output_feedback_control_with_local_measurement}, one can define the triple $(M_{1i},M_{2i},M_{3i})$ in \eqref{state_equation_of_dynamic_output_feedback_control_with_local_measurement} as 
\begin{eqnarray}
M_{1i} &\triangleq&  \begin{bmatrix}
A_i+B_i K_{1i}-H_i(C_{\rom{m}i}+D_{\rom{m}i}K_{1i}) \  & (B_i-H_iD_{\rom{m}i})K_{2i} \\
0  & G_{1i}
\end{bmatrix}, \nonumber \\
M_{2i} &\triangleq& \begin{bmatrix}
0\\
G_{2i}  
\end{bmatrix}, \ 
M_{3i} \triangleq \begin{bmatrix}
H_i\\
0  
\end{bmatrix}. \label{M1i_M2i_M3i_definition_for_the_second_controller}
\end{eqnarray}	
Using \eqref{measured_output_1} and \eqref{output_equation_of_dynamic_output_feedback_control_with_local_measurement_specified}, \eqref{local_state_observer_1} can be rewritten as
\begin{eqnarray}
\dot {\hat x}_i(t)  &=& H_iC_{\rom{m}i} x_i(t) + (A_i+B_i K_{1i}-H_iC_{\rom{m}i}){\hat x}_i(t)  +  B_iK_{2i} \bar z_i(t) ,   \quad {\hat x}_i(0)={\hat x}_{i0},  \quad t \geq 0.  \label{local_state_observer_3}
\end{eqnarray} 

Next, define $\hat x(t) \triangleq [\hat x_1^\rom{T}(t), \ldots, \hat x_N^\rom{T}(t)]\mT$, $\bar z(t) \triangleq [\bar z_1^\rom{T}(t), \ldots, \bar z_N^\rom{T}(t)]\mT$, and $ H \triangleq \rom{diag} (H_1, \ldots, H_N ) $.  Inserting \eqref{output_equation_of_dynamic_output_feedback_control_with_local_measurement_specified} into \eqref{local_state_2} and \eqref{local_tracking_error_1}, using \eqref{local_state_observer_3}, \eqref{partial_state_equation_of_dynamic_output_feedback_control_with_local_measurement}, and the above definitions, \eqref{local_state_2}, \eqref{state_equation_of_dynamic_output_feedback_control_with_local_measurement}, and \eqref{local_tracking_error_1} can be compactly written as 
\begin{eqnarray}
\dot x(t) &=& Ax(t)+ BK_1 \hat x(t) + BK_2 \bar z(t) + E{\omega_{{\rom{a}}}}(t), \hspace{1.29cm} \quad x(0)=x_{0}, \quad t \geq 0, \label{global_state_2}\\
\dot {\hat x}(t) &=& HC_{\rom{m}}x(t)+ (A +BK_1 -HC_{\rom{m}})\hat x(t) + BK_2 \bar z(t),  \quad \hat x(0)= \hat x_{0},  \quad t \geq 0, \label{global_observer_state_2}\\
\dot {\bar z}(t) &=& G_1\bar z(t)+G_2 e_{\rom{v}}(t),  \hspace{4.68cm} \quad \bar z(0)=\bar z_0, \quad t \geq 0, \label{global_partial_controller_state_2} \\
e(t)&=& Cx(t) + DK_1 \hat x(t) + DK_2 \bar z(t) - R_{\rom{a}}{\omega_{{\rom{a}}}}(t). \label{global_tracking_error_2}
\end{eqnarray} 
Now, insert \eqref{global_tracking_error_2} into \eqref{e_v} and replace the obtained expression with the one in \eqref{global_partial_controller_state_2}. Let $\eta(t) \triangleq [x^\rom{T}(t), \hat x^\rom{T}(t), \bar z^\rom{T}(t)]\mT \in \IR^{\bar n +\bar n_{z_{2}}}$, where $\bar n_{z_{2}}= \sum_{i=1}^{N} n_{z_{2i}}$. Then, the closed-loop system of \eqref{local_state_2}-\eqref{virtual_local_error_1} and  \eqref{measured_output_1}-\eqref{state_equation_of_dynamic_output_feedback_control_with_local_measurement} can be represented as 
\begin{eqnarray}
\dot \eta(t) &=& A_{\eta} \eta(t)+B_{\eta} {\omega_{{\rom{a}}}}(t),  \quad \eta(0)=\eta_0, \quad t  \geq 0, \label{output_feedback_local_measuremnt_global_state}\\
e(t)&=& C_{\eta}\eta(t)+ D_{\eta} {\omega_{{\rom{a}}}}(t), \label{output_feedback_local_measuremnt_global_error}
\end{eqnarray}
where 
\begin{eqnarray}
A_{\eta} \hspace{-0.02 cm} &=& \hspace{-0.05 cm} \begin{bmatrix}
A  & BK_1 & BK_2 \\
HC_{\rom{m}}  & A +BK_1 -HC_{\rom{m}} & BK_2 \\
G_2 \mathcal{W}C & G_2\mathcal{W}DK_1 & G_1 +G_2\mathcal{W}DK_2
\end{bmatrix}, \nonumber \\ 
B_{\eta} &=&  \begin{bmatrix}
E   \\
0 \\
-G_2\mathcal{W} R_{\rom{a}} 
\end{bmatrix},  \
C_{\eta} = \begin{bmatrix}
C & DK_1 & DK_2  
\end{bmatrix}, \ D_{\eta} = - R_{\rom{a}}. \nonumber 
\end{eqnarray}

For the following result, we define $A_{\rom{H}i} \triangleq A_i - H_i C_{\rom{m}i}$ and $A_{\rom{H}} \triangleq A-H C_{\rom{m}} $. By Assumption \ref{A7}, $H_i$ can always be chosen such that $A_{\rom{H}i} $ is Hurwitz for all $i \in \mathcal{N}$.

\theorem\label{Theorem3} Let Assumptions \ref{A1}-\ref{A3} and \ref{A6} hold. If $A_{\rom{g}}$ is Hurwitz and $A_{\rom{H}i}$ is Hurwitz for all $i \in \mathcal{N}$, then the distributed dynamic output feedback control with local measurement given by \eqref{output_equation_of_dynamic_output_feedback_control_with_local_measurement} and \eqref{state_equation_of_dynamic_output_feedback_control_with_local_measurement} solves the problem in Definition \ref{Def2}.

\textit{Proof.} Let $K \triangleq [K_1 \ K_2]$, $\hat A \triangleq A$, $\hat B \triangleq BK$, $\hat C \triangleq C$, $\hat C_{\rom{m}} \triangleq C_{\rom{m}} $, $\hat D \triangleq DK$, $\hat D_{\rom{m}} \triangleq D_{\rom{m}}K $,  $\hat E \triangleq E$, $\hat F \triangleq - R_{\rom{a}},$
\begin{eqnarray}
M_{1} &\triangleq&  \begin{bmatrix}
A+B K_{1}-H(C_{\rom{m}}+D_{\rom{m}}K_{1}) \  & (B-HD_{\rom{m}})K_{2} \\
0  & G_{1}
\end{bmatrix}\hspace{-0.05cm}, \nonumber \\
M_{2} &\triangleq& \begin{bmatrix}
0\\
G_{2}  
\end{bmatrix}, \ 
M_{3} \triangleq \begin{bmatrix}
H\\
0  
\end{bmatrix}. \label{M1_M2_M3_theorem3} 
\end{eqnarray}	
Now, observe that the quadruple $(A_{\eta},B_{\eta},C_{\eta},D_{\eta})$ takes the form of $(A_{\rom{c}},B_{\rom{c}},C_{\rom{c}},D_{\rom{c}})$ in Lemma \ref{Lemma3}. Recall from the proof of Theorem \ref{Theorem1} that the triple $(G_1,G_2,0)$ incorporates an $Np$-copy internal model of $A_{0\rom{a}}$ under Assumption \ref{A6}. This clearly implies that the triple $(M_1,M_2,M_3)$ also incorporates an $Np$-copy internal model of $A_{0\rom{a}}$. It is given that Assumptions \ref{A1} and \ref{A3} hold. 
In order to apply Lemma \ref{Lemma3}, we need to show that $A_{\eta}$ is Hurwitz under the conditions that $A_{\rom{g}}$ is Hurwitz and $A_{\rom{H}i}$ is Hurwitz for all $i \in \mathcal{N}$. To this end, the following elementary row and column operations are performed on $A_{\eta}$. First, subtract row 1 from row 2 and add column 2 to column 1. Second, interchange rows 2 and 3, and interchange columns 2 and 3. Thus, we obtain the matrix given by
\begin{eqnarray}
\bar A_{\eta} \triangleq \begin{bmatrix}
A+BK_1   & BK_2  & BK_1 \\
G_2 \mathcal{W}(C+DK_1)  & G_1 +G_2\mathcal{W}DK_2 & G_2\mathcal{W}DK_1 \\
0 & 0 & A_{\rom{H}}
\end{bmatrix}\hspace{-0.03cm}.\label{similar_to_the_original_one} \nonumber
\end{eqnarray}
Considering the performed elementary row and column operations, one can verify that $A_{\eta}$ is similar to $\bar A_{\eta}$; hence, they have the same eigenvalues. Since $\bar A_{\eta}$ is upper block triangular, 
$\sigma(\bar A_{\eta}) = \sigma(A_{\rom{g}}) \cup \sigma(A_{\rom{H}}) $. Note that $A_{\rom{H}}$ is Hurwitz as $A_{\rom{H}i}$ is Hurwitz for all $i \in \mathcal{N}$. It is also given that $A_{\rom{g}}$ is Hurwitz. Thus, $A_{\eta}$ is Hurwitz. Then, the matrix equations
\begin{eqnarray}
X_{\eta}A_{0\rom{a}} &=& A_{\eta}X_{\eta} + B_{\eta}, \label{lme1_output_feedback_local_measurement}   \nonumber  \\
0 &=& C_{\eta}X_{\eta} + D_{\eta}, \label{lme2_output_feedback_local_mesurement} \nonumber 
\end{eqnarray}
have a unique solution $X_{\eta}$ by Lemma \ref{Lemma3}.

Following similar steps to those in the proof of Theorem \ref{Theorem1}, it can be shown under Assumption \ref{A2} that $e_i(t)$ is ultimately bounded with an ultimate bound for all $\eta_0 $ and for all $i \in \mathcal{N}$. If, in addition, $\lim_{t\to\infty}  A_0 \omega(t) - \dot \omega(t)= 0$, then for all $ \eta_{0} \ $ $\lim_{t\to\infty} e_i(t)= 0, \ \forall i \in \mathcal{N} $. \hfill$\blacksquare$ 

\remark\label{Remark6} Since the condition on $A_{\rom{H}i}$ is both agent-wise and graph-wise local, obtaining an agent-wise local sufficient condition that ensures the property a) of Definition \ref{Def2} boils down to finding an agent-wise local sufficient condition, under standard assumptions, for the stability of $A_{\rom{g}}$, which is already given in Theorem \ref{Theorem2}.

\subsection{Dynamic Output Feedback}\label{third_controller}

Define $z_i(t)$, $\bar K_i$, and $u_i(t)$ as in Section \ref{second_controller}; that is, \eqref{output_equation_of_dynamic_output_feedback_control} has the form \eqref{output_equation_of_dynamic_output_feedback_control_with_local_measurement_specified}. 
Since $e_{\rom{v}i}(t)$ is the only available information to each agent,  the following distributed observer is considered instead of \eqref{local_state_observer_2} to estimate the state $x_i(t)$
\begin{eqnarray}
\dot {\hat x}_i(t)  &=& \big(A_i+B_i K_{1i}-L_i(C_i+D_iK_{1i})\big){\hat x}_i(t) +L_ie_{\rom{v}i}(t) +  (B_i-L_iD_i)K_{2i} \bar z_i(t), \nonumber \\ && \hspace{7.7 cm} \quad {\hat x}_i(0)={\hat x}_{i0},  \quad t \geq 0,  \label{state_observer_1}
\end{eqnarray} 
where $L_i$ is the observer gain matrix. Let $\bar z_i(t)$ satisfy the dynamics in \eqref{partial_state_equation_of_dynamic_output_feedback_control_with_local_measurement}. We can now define the pair $(M_{1i},M_{2i})$ in \eqref{state_equation_of_dynamic_output_feedback_control} by replacing the triple $(H_i,C_{\rom{m}i},D_{\rom{m}i})$ in $M_{1i}$ (respectively, the zero matrix in $M_{2i}$) given by \eqref{M1i_M2i_M3i_definition_for_the_second_controller} with $(L_i,C_i,D_i)$ (respectively, $L_i$).

Define $\hat x(t)$ and $\bar z(t)$ as in the previous subsection and   $ L \triangleq \rom{diag} (L_1, \ldots, L_N ).$ 
Inserting \eqref{output_equation_of_dynamic_output_feedback_control_with_local_measurement_specified} into \eqref{local_state_2} and \eqref{local_tracking_error_1}, using \eqref{state_observer_1}, \eqref{partial_state_equation_of_dynamic_output_feedback_control_with_local_measurement}, and the above definitions, \eqref{local_state_2}, \eqref{state_equation_of_dynamic_output_feedback_control}, and \eqref{local_tracking_error_1} can be expressed by \eqref{global_state_2}, 
\begin{eqnarray}
\dot {\hat x}(t) &=&  \big(A+B K_{1}-L(C+DK_{1})\big){\hat x}(t) + (B-LD)K_{2} \bar z(t) +Le_{\rom{v}}(t),  \quad \hat x(0)= \hat x_{0},  \quad t \geq 0, \label{global_observer_state_3}
\end{eqnarray} 
\eqref{global_partial_controller_state_2}, and  \eqref{global_tracking_error_2}. Next, insert \eqref{global_tracking_error_2} into \eqref{e_v} and replace the obtained expression not only with the one in \eqref{global_partial_controller_state_2} but also with the one in \eqref{global_observer_state_3}. In addition, define $\eta(t)$ as in Section \ref{second_controller}. Then, the closed-loop system of \eqref{local_state_2}-\eqref{virtual_local_error_1},  \eqref{output_equation_of_dynamic_output_feedback_control}, and \eqref{state_equation_of_dynamic_output_feedback_control} can be expressed by \eqref{output_feedback_local_measuremnt_global_state} and \eqref{output_feedback_local_measuremnt_global_error} if the second row of $A_{\eta}$ is replaced with 
\begin{eqnarray}
\begin{bmatrix}
L\mathcal{W}C \  & A +BK_1 - L(C+DK_1- \mathcal{W}DK_1) \ & (B-LD+L\mathcal{W}D)K_2 \nonumber 
\end{bmatrix}
\end{eqnarray}
and the second row of $B_{\eta}$ is replaced with $-L\mathcal{W} R_{\rom{a}} $.

\theorem\label{Theorem4} Let Assumptions \ref{A1}-\ref{A3} and \ref{A6} hold. If the resulting $A_{\eta}$ is Hurwitz, then the distributed dynamic output feedback control given by \eqref{output_equation_of_dynamic_output_feedback_control} and \eqref{state_equation_of_dynamic_output_feedback_control} solves the problem in Definition \ref{Def2}.

\textit{Proof.} Define $K$, $\hat A$, $\hat B$, $\hat C$, $\hat D$, $\hat E$, and $\hat F$ as in the proof of Theorem \ref{Theorem3}. Let $\hat C_{\rom{m}} \triangleq 0 $, $\hat D_{\rom{m}} \triangleq 0 $, and $M_3 \triangleq 0$. Define also the pair  $(M_1,M_2)$ by replacing the triple $(H,C_{\rom{m}},D_{\rom{m}})$ in $M_1$  (respectively, the zero matrix in $M_{2}$) given by \eqref{M1_M2_M3_theorem3} with $(L,C,D)$ (respectively, $L$). Then, observe that the resulting quadruple $(A_{\eta},B_{\eta},C_{\eta},D_{\eta})$ takes the form of $(A_{\rom{c}},B_{\rom{c}},C_{\rom{c}},D_{\rom{c}})$ in Lemma \ref{Lemma3}. By the same argument in the proof of Theorem \ref{Theorem3}, the resulting triple $(M_1,M_2,M_3)$ incorporates an $Np$-copy internal model of $A_{0\rom{a}}$ under Assumption \ref{A6}. Since, in addition, Assumptions \ref{A1}-\ref{A3} hold and $A_{\eta}$ is Hurwitz, the rest of the proof can be completed by following the steps given in the proof of Theorem \ref{Theorem1}. \hfill$\blacksquare$ 

Now, our goal is to obtain an agent-wise local sufficient condition that assures the property $a)$ of Definition \ref{Def2} under some standard assumptions. For this purpose, define $\mu_i(t)$ as in Section \ref{first_controller} and  let $\zeta_i (t)\triangleq [x_i^\rom{T}(t), \hat x_i^\rom{T}(t), \bar z_i^\rom{T}(t)]^\rom{T} \in \IR^{n_i +n_{z_{2i}}} $, 
\begin{eqnarray}
A_{\rom{F}i} &\triangleq& \ \begin{bmatrix}
A_i  & B_iK_{1i} & B_iK_{2i} \\
L_iC_i  & A_i +B_iK_{1i} -L_iC_i & B_iK_{2i} \\
G_{2i}C_i & G_{2i}D_iK_{1i} & G_{1i} +G_{2i}D_iK_{2i}
\end{bmatrix}, \nonumber \\
B_{\rom{F}i} &\triangleq& \ \begin{bmatrix}
0 \\
-L_i \\
-G_{2i} 
\end{bmatrix},  \
C_{\rom{F}i} \triangleq \begin{bmatrix}
C_i & D_iK_{1i} & D_iK_{2i}  
\end{bmatrix}. \nonumber 
\end{eqnarray}
Furthermore, consider \eqref{local_state_2}, \eqref{state_equation_of_dynamic_output_feedback_control}, \eqref{virtual_local_error_2}, and \eqref{local_tracking_error_1}  when $\omega(t)  \equiv 0$. By inserting \eqref{output_equation_of_dynamic_output_feedback_control} into the considered equations, we have
\begin{eqnarray}
\dot \zeta_i(t) &=& A_{\rom{F}i} \zeta_i(t) + B_{\rom{F}i}  \mu_i(t),  \quad \zeta_i(0)=\zeta_{i0}, \quad t \geq 0, \label{local_closed_loop_state1} \\
e_i(t)&=& C_{\rom{F}i} \zeta_i(t). \label{local_closed_loop_output2}
\end{eqnarray}

\remark\label{Remark7} Let $A_{\rom{L}i} \triangleq A_i - L_i C_{i}$. By performing the elementary row and column operations given in the proof of Theorem \ref{Theorem3} on $A_{\rom{F}i}$, one can show that  $\sigma( A_{\rom{F}i}) = \sigma(A_{\rom{f}i}) \cup \sigma(A_{\rom{L}i}) $.
Note that by Assumption \ref{A8}, $L_i$ can always be chosen such that  $A_{\rom{L}i}$ is Hurwitz for all $i \in \mathcal{N}$.
In conjunction with Remark \ref{Remark4}, this shows that under Assumptions \ref{A4}-\ref{A6} and Assumption \ref{A8}, it is always possible to find $K_{1i}$, $K_{2i}$, and $L_i$ such that $A_{\rom{F}i}$ is Hurwitz for all $i \in \mathcal{N}$.

Let $g_{\rom{F}i}(s) \triangleq C_{\rom{F}i}(sI -A_{\rom{F}i})^{-1}B_{\rom{F}i}$. For the dynamic output feedback case, we now state the following theorem.

\theorem\label{Theorem5} Let Assumption \ref{A3} hold and $A_{\rom{F}i}$ be Hurwitz for all $i \in\mathcal{N}$. If 
\begin{eqnarray}
\norm {g_{\rom{F}i}}_\infty \ \hspace{-0.05cm} \rho(\mathcal{F}\mathcal{A}) < 1, \quad \forall i \in \mathcal{N}, \label{local_sufficient_cond2}
\end{eqnarray}
then the resulting $A_{\eta}$ is Hurwitz.

\textit{Proof.} It follows from  Section \ref{first_controller} by comparing \eqref{local_closed_loop_state1} and \eqref{local_closed_loop_output2} with \eqref{local_closed_loop_state} and \eqref{local_closed_loop_output}. \hfill$\blacksquare$


\section{Illustrative Numerical Examples}\label{INE}

\begin{figure}[ht!] 
	\begin{center}\vspace{-0.3cm}
		\includegraphics[scale= 0.4]{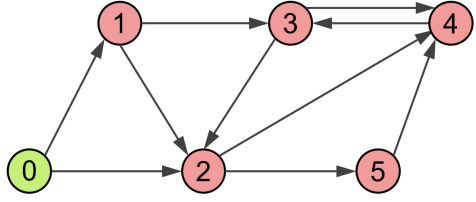}  
	\end{center} \vspace{-0.7cm}
	\caption{Augmented directed graph $\mathcal{\bar G}$. } 
	\label{Graph}
\end{figure}

To illustrate some results from the previous section, we provide two numerical examples with different exosystems. In particular, the first (respectively, second) example presents the distributed dynamic state (respectively, output) feedback control law. 
For both examples, we consider five agents with the following system, input, output, and direct feedthrough matrices 
\begin{eqnarray} 
A_i &=&
\begin{bmatrix}
-1  & 1\\
0.2  & 0 
\end{bmatrix}, \  B_i = 
\begin{bmatrix}
1    \\
2   
\end{bmatrix},\  C_i = 
\begin{bmatrix}
1  & 0  
\end{bmatrix}, \  D_i = 0.1, \nonumber \ i =1,4,5, \\
A_i &=&
\begin{bmatrix}
0  & 1 & 0\\
0  & 2 & 1 \\
0 & 0 & 0 
\end{bmatrix}, \   B_i = 
\begin{bmatrix}
0 & 0    \\
1 & 0   \\
0 & 1  
\end{bmatrix}, \   C_i = 
\begin{bmatrix}
1  & 0  & 0.4
\end{bmatrix}, \ D_i =0,  \nonumber \ i =2,3,   \nonumber 
\end{eqnarray}   
and the augmented graph $\mathcal{\bar G}$ shown in Figure \ref{Graph}. With this setup, each agent satisfies Assumptions \ref{A4} and \ref{A8}. It is also clear from Figure \ref{Graph} that Assumption \ref{A3} holds. In the simulations, we set each nonzero $a_{ij}$ to 1 and $k_i = 1, \ \hspace{-0.09cm} i=1,2$. Moreover, initial conditions for the agents are given by $x_{10} = [1,\ 0.6]^\rom{T}$, $x_{20} = [-0.5, \ 0, \ -0.2]^\rom{T}$,  $x_{30} = [-0.2, \ -0.3, \ 0]^\rom{T}$, $x_{40} = [0.6,\ 0]^\rom{T}$, $x_{50} = [0, \ 0.5]^\rom{T}$ and the controller states of all agents are initialized at zero.

\subsection{Example 1}

In this example, the disturbance $\delta(t)$ and the trajectory of the leader $r_0(t)$ satisfy the following dynamics 
\begin{eqnarray}
\dot \delta(t) &=& \begin{bmatrix}
0 & 0.01 & 0 \\
0 &  0  &  0\\
0 &  0  &  -0.05
\end{bmatrix}\delta(t)  +   \begin{bmatrix}
0  \\
0  \\
0.05   
\end{bmatrix}, \quad \delta(0)  =   \begin{bmatrix}
0 \\
-0.2 \\
0   
\end{bmatrix}, \quad  t\geq 0,   \nonumber \\
\dot r_0(t) &=& -r_0^3(t) + u_0(t), \quad r_0(0)=0, \quad t\geq 0,  \nonumber
\end{eqnarray}
respectively, where
\begin{subnumcases}{u_0(t) = }
0.1t, \hspace{4.85cm}  0 \leq t <100, \nonumber
\\
0.1t - 2\rom{sin}(0.1t)e^{-0.01(t-100)}, \quad  100 \leq t <200, \nonumber
\\
14 + \rom{sin}(0.05(t - 200)), \hspace{2.4cm} t \geq 200 \nonumber.
\end{subnumcases} 
By the solution of the disturbance dynamics with the given initial condition, $\dot \delta(t)$ is bounded. Since $u_0(t)$ is piecewise continuous and bounded, $r_0(t)$ is bounded by Example 4.25 in \cite{Khalil_book}; hence, $\dot r_0(t)$ is piecewise continuous and bounded. Clearly, $\dot \omega(t)$ is piecewise continuous and  bounded.  
Furthermore, the exosystem affects the state of each agent and its tracking error through matrices
\begin{eqnarray} 
E_{\delta_1}   &=&
\begin{bmatrix}
0 & 1 & 0\\
0 &  0 & 0
\end{bmatrix}\, \  E_{\delta_4}   =
\begin{bmatrix}
0.1 & 0 & 0\\
0 &  0 & -0.1
\end{bmatrix}, \  E_{\delta_5}   =
\begin{bmatrix}
0 & 0 & 0\\
-0.1 &  -0.2 & 0
\end{bmatrix}, \nonumber \\
E_{\delta_2}  &=& \begin{bmatrix}
0 & 0 & 1\\
0 &  0 & 0 \\
0 & 0 & 0.5
\end{bmatrix}, \   E_{\delta_3}  = \begin{bmatrix}
0 & -0.5 & 0\\
0 &  0 & -1 \\
0 & 0.4 & 0
\end{bmatrix},  \ \nonumber R_{\rom{r}} = 1.
\end{eqnarray}  

\begin{figure}[ht!] 
	\centering
	\includegraphics[scale= 0.75]{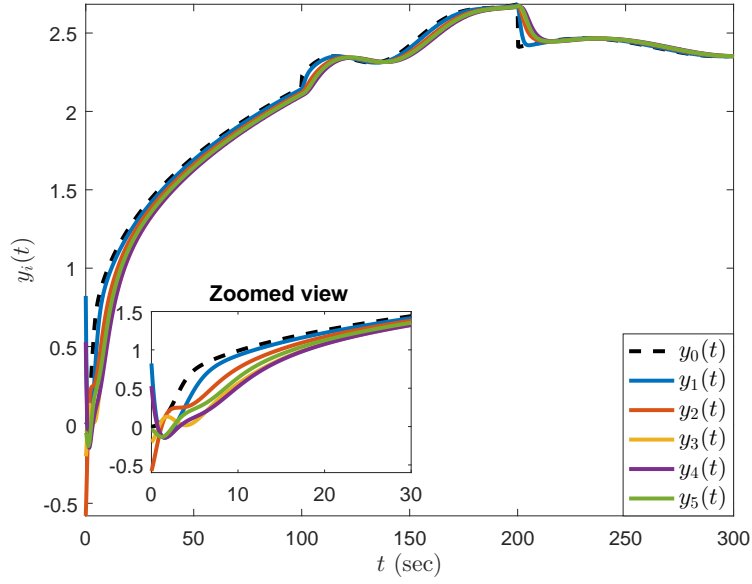} \hspace{0.0cm}
	\vspace{-0.4cm} \caption{Output responses of the agents in Example 1. } 
	\label{Example1}
\end{figure}

Suppose the piecewise continuity and boundedness of $\dot \omega(t)$ are the only information that we know about the exosystem. As it is suggested in the part $a)$ of Remark \ref{Remark2}, we then let $A_0=0$ and $(G_{1i},G_{2i})= (0,1)$ for all $i \in \mathcal{N}$. Thus, Assumptions \ref{A1}, \ref{A2}, \ref{A5}, and \ref{A6} hold. With the following controller parameters 
\begin{eqnarray} 
K_{1i} &=&-
\begin{bmatrix}
1.1960  & 0.9611
\end{bmatrix}, \ K_{2i} = -1.4142, \nonumber  \ i = 1,4,5,\\
K_{1i} &=&-
\begin{bmatrix}
4.2328  & 5.3904 &1.4038 \\
1.2604  &  1.4038 &  1.7115
\end{bmatrix}, \ K_{2i} =- \begin{bmatrix}
1.2788 \\
1.3655
\end{bmatrix},  \nonumber \ i =2,3,
\end{eqnarray}  
$A_{\rom{f}i}$ is Hurwitz for all $i \in \mathcal{N}$ and the condition given by \eqref{local_sufficient_cond} is satisfied. Thus, $A_{\rom{g}}$ is Hurwitz by Theorem \ref{Theorem2}. As Theorem \ref{Theorem1} promises, ultimately bounded tracking error is observed in Figure \ref{Example1}.

\subsection{Example 2}

The disturbance and the trajectory of the leader satisfy
\begin{eqnarray}
\dot \delta(t) &=& e^{-0.1t}, \quad \delta(0)=1, \quad t\geq 0,  \nonumber \\
\dot r_0(t) &=& \begin{bmatrix}
0 & 0.5  \\
-0.5 &  0  
\end{bmatrix}r_0(t)  +  \begin{bmatrix}
t e^{-t} \rom{sin}(t)\\
2 e^{-t}
\end{bmatrix}, \ r_0(0)  =   \begin{bmatrix}
-1 \\
1 
\end{bmatrix}, \  t\geq 0,  \nonumber
\end{eqnarray}
respectively. Moreover, $E_{\delta_1} = [1 \ 0]\mT$, $E_{\delta_2} = [0 \ 1 \ 0]\mT$, $E_{\delta_3} = [-1.5 \ 0 \ 0.3]\mT$, $E_{\delta_4} = [0 \ 2]\mT$, \\ $E_{\delta_5} = [0.2 \ -0.2]\mT$, and $R_{\rom{r}}= [1 \ 0]$.

\begin{figure}[hb!] 
	\centering
	\includegraphics[scale= 0.75]{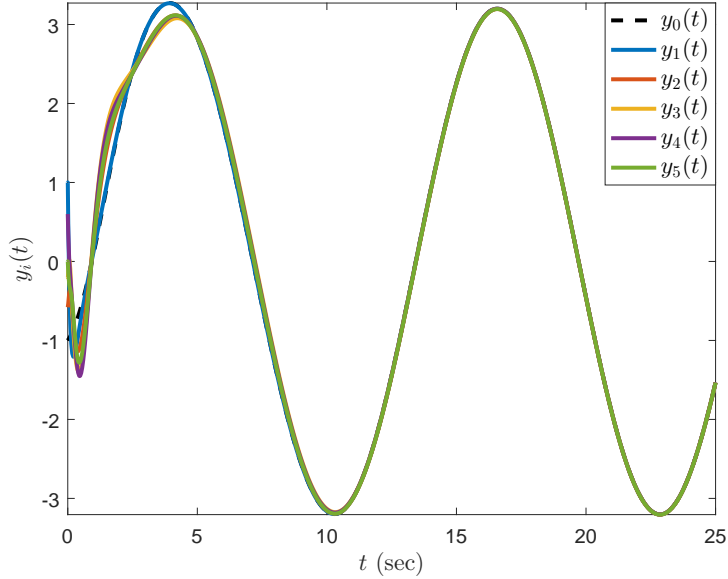} \hspace{0.0cm}
	\vspace{-0.4cm} \caption{Output responses of the agents in Example 2. } 
	\label{Example2}
\end{figure}

Suppose the unforced parts of the given dynamics are available to a control designer and the forcing terms are known to be piecewise continuous and convergent to zero. Then, let 
\begin{eqnarray} 
A_0 &=&
\begin{bmatrix}
0  & 0.5 & 0 \\
-0.5  &0  & 0  \\
0  & 0 & 0
\end{bmatrix}, \nonumber 
\end{eqnarray}  
and
\begin{eqnarray} 
G_{1i} &=&
\begin{bmatrix}
0  & 1 & 0 \\
0  &0  & 1  \\
0  & -0.25 & 0
\end{bmatrix}, \ G_{2i} =  \begin{bmatrix}
0   \\
0  \\
1
\end{bmatrix}, \nonumber \quad \forall i \in \mathcal{N}.
\end{eqnarray}  
Hence, Assumptions \ref{A1}, \ref{A5}, and \ref{A6} hold. In addition, $\lim_{t\to\infty}  A_0 \omega(t) - \dot \omega(t)= 0$. 
Note that Assumption \ref{A2} automatically holds since $A_0 \omega(t) - \dot \omega(t)$ is piecewise continuous and convergent. With the following controller parameters 
\begin{eqnarray} 
K_{1i} &=&-
\begin{bmatrix}
5.1794  & 0.7932
\end{bmatrix}, \ L_{i} = \begin{bmatrix}
17  & 80.2
\end{bmatrix}^{\rom{T}},  \nonumber \\
K_{2i} &=&- 
\begin{bmatrix}
2 & 5.4458 &  10.3182
\end{bmatrix},   \  i = 1,4,5, \nonumber \\
K_{1i} &=&-
\begin{bmatrix}
6.1916   & 5.7686   &1.7835\\
3.9299  &  1.7835 &  2.4282
\end{bmatrix}, \  L_{i} = \begin{bmatrix}
-187  & 756 & 600
\end{bmatrix}^{\rom{T}},     \nonumber  \\
K_{2i} &=&-
\begin{bmatrix}
0.4513  & 0.9173 &  3.3839\\
0.8924  & 2.2285   & 5.6377
\end{bmatrix},  \nonumber \ i =2,3,
\end{eqnarray}  
$A_{\rom{F}i}$ is Hurwitz for all $i \in \mathcal{N}$ and the condition given by \eqref{local_sufficient_cond2} is satisfied. Thus, $A_{\eta}$ is Hurwitz by Theorem \ref{Theorem5}. Furthermore, it is guaranteed by Theorem \ref{Theorem4} that $\lim_{t\to\infty} e_i(t) = 0, \ \forall i\in\mathcal{N}$ and this fact is demonstrated in Figure \ref{Example2}.

\section{Conclusion}\label{con}

In this paper, we studied the cooperative output regulation problem of heterogeneous linear time-invariant multiagent systems over fixed directed communication graph topologies. Specifically, we introduced a new definition of the linear cooperative output regulation problem (see Definition \ref{Def2}), which allows a broad class of functions to be tracked and rejected by a network of agents, and focused on an internal model based distributed control approach. For the three different distributed control laws (i.e., dynamic state feedback, dynamic output feedback with local measurement, and dynamic output feedback), we investigated the solvability of this problem, which resulted in global and local sufficient conditions (see Theorems \ref{Theorem1}-\ref{Theorem5}).  In addition, the provided two numerical examples illustrated the efficacy of our contributions. Finally, we reported and addressed a considerable number of gaps in the existing related literature (see Appendices and Section \ref{first_controller}).


\section*{Acknowledgment}

The authors would like to thank Dr. Chao Huang for helpful responses to our questions regarding the results in \cite{Huang_C_tac}. 

\section*{Appendices }\label{Appendix} 

\subsection*{Appendix A. Solvability of \eqref{lme1_state_feedback} and \eqref{lme2_state_feedback}}

Section III in \cite{Huang_C_tac} also studies the solvability of the matrix equations in \eqref{lme1_state_feedback} and \eqref{lme2_state_feedback}, which correspond to the matrix equations given by (6) in \cite{Huang_C_tac}, with an alternative approach.
Specifically, the last paragraph of Section III in \cite{Huang_C_tac} lists three sufficient conditions based on Remark 3.8 of \cite{Huang_TAC_1994} to guarantee that these matrix equations have a unique solution. However, it cannot be guaranteed as it is claimed in \cite{Huang_C_tac}. This section aims to present the gaps between the conditions and the existence of a unique solution to the matrix equations, propose appropriate modifications that fill these gaps, and explain the motivation behind our approach.  
For this purpose, we first focus on Definition 3.7 and Remark 3.8 in \cite{Huang_TAC_1994} to fix a problem in \cite{Huang_TAC_1994}. Then, we revisit the conditions listed in \cite{Huang_C_tac} to point out the missing one.
Finally, a motivational example is provided and the difference between the approach in \cite{Huang_C_tac} and the one in this paper is highlighted. 

\vspace{0.2cm}

\textit{In this paragraph, the notation and the terminology in \cite{Huang_TAC_1994} are adopted and readers are referred to (3.5), (3.6), (3.8), Definition 3.7, and Remark 3.8 in \cite{Huang_TAC_1994}.} The problem in \cite{Huang_TAC_1994} is that the conditions of Remark 3.8 do not ensure the stabilizability of the pair given by (3.8). Moreover, this problem is directly transferred to \cite{Huang_C_tac}. To illustrate this point,  we consider the following system, input, output, and direct feedthrough matrices of the plant; and system matrix of the exosystem
\begin{eqnarray} 
A =
\begin{bmatrix}
1  & 2\\
1  & 0 
\end{bmatrix}, \ B = 
\begin{bmatrix}
2   \\
0     

\end{bmatrix},\ C = 
\begin{bmatrix}
0.5  & -0.5  
\end{bmatrix}, \ D = 0, \ A_1 = 0. \nonumber
\end{eqnarray}   
It can be easily checked that the plant and the exosystem above satisfy the first and the second conditions of Remark 3.8. Note that $m(s)=s$ is the minimal polynomial of $A_1$. Then, choose the pair $(\beta_1,\sigma_1)$ in (3.6) as follows 
\begin{eqnarray} \beta_1 =
\begin{bmatrix}
0  & 1\\
0  & 1 \\
\end{bmatrix}, \ \sigma_1 = 
\begin{bmatrix}
0 \\
1  \\
\end{bmatrix}. \nonumber
\end{eqnarray}
It is obvious that the pair $(\beta_1,\sigma_1)$ is controllable and the minimal polynomial of $A_1$ divides the characteristic polynomial of $\beta_1$. Thus, the pair $(G_1,G_2) \triangleq (\beta_1,\sigma_1)$ incorporates a 1-copy internal model of $A_1$ according to Definition 3.7. Let us now investigate the stabilizability of the pair in (3.8). This pair is not controllable by the controllability matrix test (e.g., see Theorem 12.1 in \cite{Hespanha_book}) and the eigenvalues of the first matrix of this pair are $-1$, $0$, $1$, and $2$. The eigenvector test for stabilizability (e.g., see Theorem 14.1 in \cite{Hespanha_book}) reveals that unstable eigenvalue $1$ is the uncontrollable mode; that is, the pair in (3.8) is not stabilizable. Hence, there do not exist $K_1$ and  $K_2$ such that $A_c$ defined in (3.5) is Hurwitz. This counterexample to Remark 3.8 is obtained due to the fact that the constructed $G_1$ violates Property 1.5 in \cite{Huang_book}. In fact, J. Huang (personal communication, June 9, 2018) recognizes the problem in Remark 3.8; hence, he adds Property 1.5 as a condition to Lemma 1.26\footnote{\hspace{0.16cm}We also note that the proof of Lemma 1.26 in \cite{Huang_book} is still valid even if Assumption 1.1 in \cite{Huang_book} is removed from the hypotheses of Lemma 1.26.} of \cite{Huang_book}.

In this paper, Definition \ref{Def1} modifies the second property of Definition 1.22 given after (1.58) in \cite{Huang_book}. This modification guarantees that Property 1.5 in \cite{Huang_book} automatically holds if Assumption \ref{A5} holds. Based on the foregoing discussions, it is clear that Remark \ref{Remark4} is true. 

\textit{The following two paragraphs adopt the notation and the terminology from \cite{Huang_C_tac}. Readers are referred to (5), (6), (7), (8), (10), Definition 2, Lemma 2, Section II.B, and Section III in \cite{Huang_C_tac}.}
It is shown in Section III that if the matrix equations in (8) have solutions $X_{1i}$ and $X_{2i}$ for $i = 1, \ldots, N$, then the ones in (7) have  solutions $X_1 = \rom{diag} (X_{11}, \ldots, X_{1N} ) $ and $X_2 = \rom{diag} (X_{21}, \ldots, X_{2N} ) $; that is, the matrix equations in (6) have a solution $X = [X_1^{\rom{T}} \ X_2^\rom{T}]\mT $. Furthermore, it is claimed that if the three conditions\footnote{\hspace{0.16cm}In Section II.B, $S$ is assumed to have no strictly stable modes.} listed in the last paragraph of Section III hold, then the matrix equations in (8) have unique solutions $X_{1i}$ and $X_{2i}$ for $i = 1, \ldots, N$. However, these conditions do not guarantee the unique solutions. For, consider
$A_1 = 0$, $B_1 = 1$, $C_1 = 1$, $D_1=0$, $S=0$, $R =1$, $P_1 = 1$, $F_1 = 0$, and $G_1 =1$. It can be easily checked that the listed conditions are satisfied and Property 1.5 in \cite{Huang_book} is not violated. Choose $K_1=0$ and $H_1=0$. From the first matrix equation in (8), we get $1=0$, which is a contradiction. 
We now point out the problem in the claim. First, observe that the matrix equations in (8) can be equivalently written as the matrix equations given by (1.70) and (1.71) in \cite{Huang_book}. Then, by Lemma 1.27 in \cite{Huang_book}, one can note that the following condition is missed in the claim: \textit{$\tilde A_i$ given after (10) is Hurwitz\footnote{\hspace{0.16cm}After the suggested modification above, $K_i$ and $H_i$ can always be chosen such that $\tilde A_i$ is Hurwitz under the listed conditions.} for $i = 1, \ldots, N$}. 
It can be shown that this condition, together with the assumption on $S$, ensures that zero matrices are the \textit{unique} solutions to the off-block-diagonal matrix equations in (7) by adding $G_c\big((C_c+D_cK_c)X_1+D_cH_cX_2 -R_c\big)$ to the left side of the second equation in (7) that gives an equivalent form of (7) and applying the first part of Proposition A.2 in \cite{Huang_book}. In conclusion, if the assumption on $S$ holds, the third condition in the list holds for $i =1, \ldots, N$, and \textit{$\tilde A_i$ is Hurwitz for $i = 1, \ldots, N$}, then the matrix equations in (6) have a \textit{unique} solution $X$.

According to Lemma 2, the problem in Definition 2 is solved if the assumption on $S$ holds, $A_l$ given after (5) is Hurwitz, and the matrix equations in (6) have a unique solution $X$. Although the approach utilized during the derivation of the listed conditions does not take into account the assumption on $A_l$, one may wonder the answer of the following question: \textit{Let the listed conditions hold and $A_l$ be Hurwitz. Then, can we conclude that $\tilde A_i$ is Hurwitz for $i = 1, \ldots, N$?} The answer is \textit{no}. That is, the missing condition cannot be satisfied by assuming that the listed conditions hold and $A_l$ is Hurwitz. To clarify this point, consider the system parameters of the agents, the system matrix of the exosystem, and the adjacency matrix of $\mathcal{G^*}$
\begin{eqnarray} 
A_1 &=&
\begin{bmatrix}
-1  & 1\\
1  & 0 
\end{bmatrix}, \ B_1 = 
\begin{bmatrix}
1 & 0.5   \\
0 & 0.25    
\end{bmatrix},\ C_1 = 
\begin{bmatrix}
1  & -0.5  
\end{bmatrix}, \ D_1 = 0, \nonumber 
\end{eqnarray}
\begin{eqnarray} 
A_2 &=&
\begin{bmatrix}
0  & 1 & 0\\
0  & 0 & 1 \\
0 & 0 & 0 
\end{bmatrix}, \ B_2 = 
\begin{bmatrix}
0   \\
0 \\ 
1
\end{bmatrix},\ C_2 = 
\begin{bmatrix}
1  & 0 & 0  
\end{bmatrix}, \ D_2 = 0, \nonumber	
\\
A_3 &=& 1, \ B_3 = -1, \ C_3 =1,\ D_3 =0,\ S=0, \nonumber		 \\
Q^* &=&
\begin{bmatrix}
1  & 0 & 0 & 0\\
0.5 & 0 & 0 &0.5 \\
0 & 0.5 & 0 & 0.5 \\
0 & 0.5 & 0.5 & 0
\end{bmatrix}. \nonumber 
\end{eqnarray}   
Choose $(F_{i},G_{i}) = (0,1)$, $i=1,2,3.$ It can be easily checked that the listed conditions are satisfied and Property 1.5 in \cite{Huang_book} is not violated. One can also obtain $W$, which is required to construct $A_l$, from $Q^*$. Then, choose the remaining parameters of the controllers as follows
\begin{eqnarray} 
K_1 &=&
\begin{bmatrix}
2.6752   &  9.6624\\
-10.6752  & -24.6624 
\end{bmatrix}, \ H_1 = 
\begin{bmatrix}
-6.4   \\
6.4    
\end{bmatrix}, \nonumber \\
K_2 &=&-
\begin{bmatrix}
104.56  & 57.936 & 14.828
\end{bmatrix},  \ H_2 = -80,  \ K_3 = 0.8,  \ H_3 =1.   \nonumber 
\end{eqnarray}   
With this setup, it can be verified that $\tilde A_3$ is not Hurwitz even though $A_l$ is Hurwitz. 

Based on the previous example, the following question arises: \textit{Is the missing condition in \cite{Huang_C_tac} necessary to ensure that the matrix equations given by (6) in \cite{Huang_C_tac} have a unique solution?} In fact, this question is the motivation behind the key lemma (i.e., Lemma \ref{Lemma3}) of this paper and the answer is \textit{no}. In contrast to Section III in \cite{Huang_C_tac}, the approach in Lemma \ref{Lemma3} does not decompose matrix equations, which consist of the overall dynamics of the multiagent system, into matrix equations, which deal with the dynamics of each agent separately; hence, the missing condition in \cite{Huang_C_tac} is not required in Lemma \ref{Lemma3}. Furthermore, 
not only dynamic state feedback but also dynamic output feedback with local measurement and dynamic output feedback effectively utilize Lemma \ref{Lemma3} to solve the stated problem in Definition \ref{Def2}  (see Theorems \ref{Theorem1}, \ref{Theorem3}, and \ref{Theorem4}).

Since the proof of Theorem 1 and the statement of Theorem 4 in \cite{Farnaz_automatica} use the approach in Section III of \cite{Huang_C_tac}, we believe that the discussion in this section will also be helpful for the readers of the results in \cite{Farnaz_automatica}.

\subsection*{Appendix B. On Theorem 2 in \cite{Huang_C_tac} }

\textit{In this section, the notation and the terminology in  \cite{Huang_C_tac} are adopted and readers are referred to (5), (10), (15), and Theorem 2 in  \cite{Huang_C_tac}.}
Now, consider the system parameters of the agent, the system matrix of the exosystem, and the adjacency matrix of $\mathcal{G^*}$ given by
\begin{eqnarray} 
A_1 &=&
\begin{bmatrix}
1  & 0 & 0 \\
0  & 1 & 0  \\
0  & 0 & -1
\end{bmatrix}, \ B_1 = I_3, \  C_1 =
\begin{bmatrix}
1  & 0 & 0 \\
0  & 1 & 0  \\
\end{bmatrix}, \  D_1 = 0, \ S=0, \
Q^*=\begin{bmatrix}
1  & 0 \\
1  & 0   
\end{bmatrix}.  \nonumber
\end{eqnarray}
Choose $(F_1,G_1)=(0,I_2)$ and 
\begin{eqnarray} 
K_1 &=& 
\begin{bmatrix}
-2  & 0 & 0 \\
0  & -2 & 0  \\
0  & 0 & 2
\end{bmatrix},  \  H_1 =
\begin{bmatrix}
-1  & 0 \\
0  & -1   \\
0 & 0
\end{bmatrix}.  \nonumber
\end{eqnarray}
Note that $W = 1$ from $Q^{*}$; hence, $A_l$ given after (5) is nothing but $\tilde A_1$ given after (10). With this setup, one can verify that $T_1(s)$ given before Theorem 2 is stable and the condition in (15) is automatically satisfied, but $A_l$ is not Hurwitz. This counterexample is obtained because the realization of $T_1(s)$ is neither stabilizable nor detectable. In fact, a loss of one of them is enough to find a counterexample.

The above setup also applies to Theorem 5 in \cite{Farnaz_automatica} since it relies on Theorem 2 and its conditions are satisfied. It should be noted that although Assumptions 1-4 in \cite{Farnaz_automatica} and Property 1.5 in \cite{Huang_book} are not listed in the hypotheses of Theorem 5 in \cite{Farnaz_automatica}, this counterexample does not violate them.

\bibliographystyle{IEEEtran} \baselineskip 16pt
\bibliography{reference}

\newcommand{\SortNoop}[1]{}
\begin{thebibliography}{10}
\providecommand{\url}[1]{#1}
\csname url@samestyle\endcsname
\providecommand{\newblock}{\relax}
\providecommand{\bibinfo}[2]{#2}
\providecommand{\BIBentrySTDinterwordspacing}{\spaceskip=0pt\relax}
\providecommand{\BIBentryALTinterwordstretchfactor}{4}
\providecommand{\BIBentryALTinterwordspacing}{\spaceskip=\fontdimen2\font plus
\BIBentryALTinterwordstretchfactor\fontdimen3\font minus
  \fontdimen4\font\relax}
\providecommand{\BIBforeignlanguage}[2]{{%
\expandafter\ifx\csname l@#1\endcsname\relax
\typeout{** WARNING: IEEEtran.bst: No hyphenation pattern has been}%
\typeout{** loaded for the language `#1'. Using the pattern for}%
\typeout{** the default language instead.}%
\else
\language=\csname l@#1\endcsname
\fi
#2}}
\providecommand{\BIBdecl}{\relax}
\BIBdecl

\bibitem{Su_2012}
Y.~Su and J.~Huang, ``Cooperative output regulation of linear multi-agent
  systems,'' \emph{IEEE Transactions on Automatic Control}, vol.~57, no.~4, pp.
  1062--1066, 2012.

\bibitem{su_huang_2012_b}
------, ``Cooperative output regulation of linear multi-agent systems by output
  feedback,'' \emph{Systems \& Control Letters}, vol.~61, no.~12, pp.
  1248--1253, 2012.

\bibitem{Huang_C_tac}
C.~Huang and X.~Ye, ``Cooperative output regulation of heterogeneous
  multi-agent systems: An $\emph{H}_{\infty}$ criterion,'' \emph{IEEE
  Transactions on Automatic Control}, vol.~59, no.~1, pp. 267--273, 2014.

\bibitem{Li_2015}
Y.~Li, X.~Wang, J.~Xiang, and W.~Wei, ``Synchronised output regulation of
  leader-following heterogeneous networked systems via error feedback,''
  \emph{International Journal of Systems Science}, vol.~47, no.~4, pp.
  755--764, 2016.

\bibitem{Farnaz_automatica}
F.~\text{Adib Yaghmaie}, F.~L. Lewis, and R.~Su, ``Output regulation of linear
  heterogeneous multi-agent systems via output and state feedback,''
  \emph{Automatica}, vol.~67, pp. 157--164, 2016.

\bibitem{Cai_2017}
H.~Cai, F.~L. Lewis, G.~Hu, and J.~Huang, ``The adaptive distributed observer
  approach to the cooperative output regulation of linear multi-agent
  systems,'' \emph{Automatica}, vol.~75, pp. 299--305, 2017.

\bibitem{Lu_Liu_TAC_2017}
M.~Lu and L.~Liu, ``Cooperative output regulation of linear multi-agent systems
  by a novel distributed dynamic compensator,'' \emph{IEEE Transactions on
  Automatic Control}, vol.~62, no.~12, pp. 6481--6488, 2017.

\bibitem{Huang_book}
J.~Huang, \emph{Nonlinear output regulation: \text{T}heory and
  applications}.\hskip 1em plus 0.5em minus 0.4em\relax SIAM, 2004.

\bibitem{sarsilmaz_yucelen_asme_dscc_2017}
S.~B. Sarsilmaz and T.~Yucelen, ``On control of heterogeneous multiagent
  systems with unknown leader dynamics,'' in \emph{ASME Dynamic Systems and
  Control Conference}, 2017.

\bibitem{sarsilmaz_yucelen_acc_2018}
------, ``On control of heterogeneous multiagent systems: \text{A} dynamic
  measurement output feedback approach,'' in \emph{American Control
  Conference}, 2018.

\bibitem{Bernstein_book}
D.~S. Bernstein, \emph{Matrix mathematics: \text{T}heory, facts, and
  formulas}.\hskip 1em plus 0.5em minus 0.4em\relax Princeton University Press,
  2009.

\bibitem{Lewis_book}
F.~L. Lewis, H.~Zhang, K.~Hengster-Movric, and A.~Das, \emph{Cooperative
  control of multi-agent systems: {O}ptimal and adaptive design
  approaches}.\hskip 1em plus 0.5em minus 0.4em\relax Springer, 2014.

\bibitem{Khalil_book}
H.~K. Khalil, \emph{Nonlinear systems}.\hskip 1em plus 0.5em minus 0.4em\relax
  Prentice Hall, 2002.

\bibitem{Shivakumar}
P.~N. Shivakumar and K.~H. Chew, ``A sufficient condition for nonvanishing of
  determinants,'' \emph{Proceedings of the American Mathematical Society},
  vol.~43, no.~1, pp. 63--66, 1974.

\bibitem{Vidyasagar_1981}
M.~Vidyasagar, \emph{Input-output analysis of large-scale interconnected
  systems: \text{D}ecomposition, well-posedness, and stability}.\hskip 1em plus
  0.5em minus 0.4em\relax Springer-Verlag, 1981.

\bibitem{Wang_Chen_Liu_2018}
J.~Wang, K.~Chen, and Q.~Liu, ``Output consensus of heterogeneous multiagent
  systems with physical and communication graphs,'' \emph{Complexity}, 2018.

\bibitem{Hespanha_book}
J.~P. Hespanha, \emph{Linear systems theory}.\hskip 1em plus 0.5em minus
  0.4em\relax Princeton University Press, 2009.

\bibitem{Lavretsky_Wise_book}
E.~Lavretsky and K.~A. Wise, \emph{Robust and adaptive control with aerospace
  applications}.\hskip 1em plus 0.5em minus 0.4em\relax Springer, 2013.

\bibitem{Zhou_Doyle_Glover_book}
K.~Zhou, J.~C. Doyle, and K.~Glover, \emph{Robust and optimal control}.\hskip
  1em plus 0.5em minus 0.4em\relax Prentice Hall, 1996.

\bibitem{Nonnegative_1994}
A.~Berman and R.~J. Plemmons, \emph{Nonnegative matrices in the mathematical
  sciences}.\hskip 1em plus 0.5em minus 0.4em\relax Siam, 1994.

\bibitem{Wang_tac_2010}
X.~Wang, Y.~Hong, J.~Huang, and \text{Z.-P.} Jiang, ``A distributed control
  approach to a robust output regulation problem for multi-agent linear
  systems,'' \emph{IEEE Transactions on Automatic Control}, vol.~55, no.~12,
  pp. 2891--2895, 2010.

\bibitem{Huang_TAC_1994}
J.~Huang and C.~Lin, ``On a robust nonlinear servomechanism problem,''
  \emph{IEEE Transactions on Automatic Control}, vol.~39, no.~7, pp.
  1510--1513, 1994.

\end{thebibliography}

\end{document}